\documentclass[12pt]{article}
\usepackage{amsfonts, amsmath, amsthm, amssymb,stmaryrd}
\usepackage[all]{xy}

\numberwithin{equation}{subsection}
\newtheorem{theorem}[equation]{Theorem}
\newtheorem{lemma}[equation]{Lemma}

\newtheorem{cor}[equation]{Corollary}
\newtheorem{prop}[equation]{Proposition}

\theoremstyle{definition}

\newtheorem{remark}[equation]{Remark}

\def\AAA{\mathbb{A}}
\def\CC{\mathbb{C}}
\def\FF{\mathbb{F}}
\def\GG{\mathbb{G}}

\def\PP{\mathbb{P}}
\def\QQ{\mathbb{Q}}
\def\RR{\mathbb{R}}
\def\ZZ{\mathbb{Z}}

\def\del{\partial}
\def\dual{\vee}

\def\be{\mathbf{e}}
\def\bv{\mathbf{v}}
\def\bw{\mathbf{w}}
\def\bx{\mathbf{x}}

\def\calD{\mathcal{D}}
\def\calE{\mathcal{E}}
\def\calF{\mathcal{F}}
\def\calL{\mathcal{L}}
\def\calO{\mathcal{O}}
\def\calR{\mathcal{R}}
\def\calS{\mathcal{S}}

\def\gotha{\mathfrak{a}}
\def\gothgl{\mathfrak{gl}}
\def\gothsl{\mathfrak{sl}}
\def\gothg{\mathfrak{g}}
\def\gothh{\mathfrak{h}}
\def\gothm{\mathfrak{m}}
\def\gothn{\mathfrak{n}}
\def\gotho{\mathfrak{o}}
\def\gothz{\mathfrak{z}}

\DeclareMathOperator{\alg}{alg}
\DeclareMathOperator{\Aut}{Aut}
\DeclareMathOperator{\coker}{coker}

\DeclareMathOperator{\Gal}{Gal}
\DeclareMathOperator{\geom}{geom}
\DeclareMathOperator{\GL}{GL}
\DeclareMathOperator{\Hom}{Hom}
\DeclareMathOperator{\im}{im}
\DeclareMathOperator{\Irr}{Irr}
\DeclareMathOperator{\inte}{int}
\DeclareMathOperator{\loc}{loc}
\DeclareMathOperator{\rank}{rank}
\DeclareMathOperator{\rig}{rig}
\DeclareMathOperator{\Res}{Res}
\DeclareMathOperator{\SL}{SL}
\DeclareMathOperator{\speci}{sp}
\DeclareMathOperator{\Spec}{Spec}
\DeclareMathOperator{\Swan}{Swan}
\DeclareMathOperator{\qu}{qu}
\DeclareMathOperator{\Trace}{Trace}

\newcounter{fixmectr}

\begin{document}

\title{Fourier transforms and $p$-adic ``Weil II''}
\author{Kiran S. Kedlaya \\ Department of Mathematics, Room 2-165 
\\ Massachusetts
Institute of Technology \\ 77 Massachusetts Avenue \\
Cambridge, MA 02139 \\
\texttt{kedlaya@math.mit.edu}}
\date{July 20, 2005}

\maketitle

\begin{abstract}
We give a transcription into rigid ($p$-adic) cohomology of 
Laumon's proof of Deligne's ``Weil II'' theorem,
using a geometric Fourier transform in the spirit of $\calD$-modules.
This yields a complete, purely $p$-adic proof of the Weil conjectures
when combined with recent results on $p$-adic differential
equations of Andr\'e, Christol, Crew, Kedlaya, Matsuda, Mebkhout,
and Tsuzuki.
\end{abstract}

\tableofcontents

\section{Introduction}

\subsection{Who needs another proof, anyway?}

It has now been over thirty years since the last of
Weil's conjectures on the numbers of points of algebraic varieties over
finite fields was established by Deligne \cite{bib:del1}, following on
the groundbreaking work of the Grothendieck school in developing the
$\ell$-adic (\'etale)
cohomology of varieties in characteristic $p \neq \ell$. However,
while the $\ell$-adic cohomology is closely linked to the ``intrinsic''
geometry of a variety (namely its unramified covers), it is rather
poorly linked to its ``extrinsic'' geometry (namely its defining
equations). In the present era, the latter is a matter of some concern:
for instance,
the need to numerically compute 
zeta functions of varieties given by explicit equations is becoming
increasingly common, e.g., in coding theory and cryptography.

A more computationally oriented point of view, based on $p$-adic
analysis, can be seen in the first proof of the rationality of the
zeta function of a variety, given by Dwork \cite{bib:dwork}.
(Indeed, recently this proof has actually been converted into an
algorithm for computing zeta functions by Lauder and Wan
\cite{bib:lw}.) The intrinsic computability in Dwork's approach
makes it desirable
to build it into a $p$-adic cohomology theory with the same level of formal
manipulability as \'etale cohomology. 

We will describe the history
of the $p$-adic cohomological program
 in the next section; for now, we simply observe that
one benchmark of progress in $p$-adic cohomology
is whether one can recover the Weil conjectures by $p$-adic methods,
without reference to \'etale cohomology. More precisely, one would
like an analogue of ``Weil II'' (i.e., the main theorem of Deligne's
\cite{bib:deligne}), which relates the action of Frobenius on the
cohomology of a ``sheaf'' (that is, a coefficient object in the
cohomology theory) with the action of Frobenius on its fibres at 
points.

The main purpose of this paper is twofold:
to establish that it is indeed possible to derive the Weil conjectures
and establish an analogue of Weil II
purely within a $p$-adic cohomological framework, and to do so in
a relatively transparent, self-contained manner. (Beware that the
adjective ``self-contained'' refers to the proof within the context of rigid
cohomology, not to this paper in isolation; in particular, we rely heavily
on \cite{bib:me8}.)
To contrast, we mention
two earlier versions of ``$p$-adic Weil II''.
A purely $p$-adic derivation of the Weil conjectures was outlined by
Faltings \cite{bib:falt}, using a relative version of crystalline
cohomology. However, fleshing out the outline seems to present a formidable
technical challenge, and to our knowledge this has not been carried out;
in any case, it does not meet the transparency criterion.
A subsequent version of Weil II, based on the technically less challenging
rigid cohomology, was given by Chiarellotto \cite{bib:chiar}. However,
it is not purely $p$-adic:
it ultimately relies on Katz and Messing's crystalline version of
the Weil conjectures \cite{bib:km}, which in turn rests on Deligne.

\subsection{The $p$-adic cohomological program}

Ever since Dwork pioneered the use of $p$-adic methods in the study
of algebraic varieties over finite fields, authors too numerous
to list have attempted to ``complete'' Dwork's ideas into a more
comprehensive cohomology theory of algebraic varieties over fields
of positive characteristic, or to put it briefly, a
$p$-adic Weil cohomology. We quickly recall some of these efforts,
ending with the recent progress that makes our present approach
to $p$-adic Weil II feasible.

A couple of natural candidates fail to be Weil cohomologies
in part because they produce spaces which are ``too small''; these include 
$p$-adic \'etale cohomology and Serre's Witt vector cohomology.
(The failures in both cases can be ``explained'' using the de Rham-Witt
complex, as in \cite{bib:illusie}.)
More successful have been crystalline cohomology and
Monsky-Washnitzer cohomology (for more on which see below), but these are
restricted to smooth proper and smooth affine varieties,
respectively. A promising attempt to globalize the Monsky-Washnitzer
construction was pursued in
a series of papers by Lubkin
\cite{bib:lub1}, \cite{bib:lub2}, \cite{bib:lub3},
\cite{bib:lub4}, \cite{bib:lub5}, but a number of results in this theory
lack adequate proofs. (We believe
some gaps and inconsistencies in Lubkin's work
are addressed by ongoing work of Borger.)
A new theory in a similar vein (which also seems
to work well for singular varieties) is
Grosse-Kl\"onne's dagger cohomology \cite{bib:gk}; this theory, and 
to a lesser extent Lubkin's construction, turn out
to be closely related to rigid cohomology (for more on which see below).
% \cite{bib:lub6}.

Crystalline cohomology was developed chiefly by Berthelot and Ogus,
following ideas of Grothendieck, into a theory that produces finite
dimensional cohomology spaces for smooth proper varieties, and yields
Poincar\'e duality, the K\"unneth formula and the Lefschetz trace formula
for Frobenius. Moreover, the theory admits nonconstant coefficient
modules ($F$-crystals) and one has finite dimensionality for cohomology
with nonconstant coefficients. However, crystalline cohomology is not
finite dimensional in general for varieties which fail to be smooth
and proper. 
It should be possible to put together a
 relative theory that can be used in the nonproper case, as outlined by
Faltings \cite{bib:falt}; however, fleshing out this proposal seems to
require mastery of a large number of technical details which (to our 
best knowledge) do not appear in the literature.

Monsky-Washnitzer (MW) cohomology was introduced by Monsky and Washnitzer
\cite{bib:mw1}, \cite{bib:mw2}, \cite{bib:mw3} as an offshoot of
Dwork's methods; its definition is restricted to smooth affine varieties,
but there it admits a Lefschetz trace formula for Frobenius proved
using $p$-adic analytic techniques.
MW cohomology has the appealing feature of being very explicitly
constructed from the defining equations of a given variety, in a manner
similar to de~Rham cohomology. (Accordingly, it too has been used recently
to give explicit algorithms for computing zeta functions; see
\cite{bib:mecount}.)
However,
for a long time a proof of finite dimensionality of cohomology was lacking
except for curves \cite{bib:monsky}, limiting the usefulness of the theory.

The crystalline and MW points of view were reconciled magnificently by
Berthelot with the construction of rigid cohomology. This theory
coincides with (the torsion-free part of)
crystalline cohomology for smooth proper
varieties and with MW cohomology for smooth affine varieties; in fact,
it appears to be ``universal'' among $p$-adic cohomology theories with
field coefficients. (As field coefficients suffice for discussion of
the Weil conjectures, we steer clear of the thorny issue of constructing
\emph{integral} $p$-adic cohomologies.)
Berthelot proved finite dimensionality of rigid cohomology
for an arbitrary smooth
variety \cite{bib:ber2}
by reducing to the crystalline case, thus proving finite dimensionality
of MW cohomology. (The latter can also be proved directly using results
on $p$-adic differential equations; see \cite{bib:meb1}.)
Berthelot also established Poincar\'e duality in rigid cohomology
\cite{bib:ber4}; these results can be used to give purely $p$-adic 
proofs of the rationality and functional equation of the zeta function,
but they are not enough to allow Deligne's results on weights to be
transposed into a $p$-adic context.

Berthelot also introduced nonconstant coefficient objects in rigid cohomology,
known as overconvergent $F$-isocrystals. Finiteness of cohomology with
nonconstant coefficients was expected to follow from a conjecture of Crew
on quasi-unipotent $p$-adic differential equations, loosely analogous to
Grothendieck's local monodromy theorem.
Proofs of Crew's conjecture have now been given by
by Andr\'e \cite{bib:andre}, Mebkhout \cite{bib:meb2},
and the author \cite{bib:me7}; as a consequence,
finiteness of rigid cohomology (with and without compact supports)
with coefficients in an overconvergent $F$-isocrystal,
plus Poincar\'e duality and the K\"unneth formula,
have been obtained by the author \cite{bib:me8}.

The proof of finiteness of rigid cohomology with coefficients
makes it now feasible to establish an analogue of Weil II for
overconvergent $F$-isocrystals.
In fact, this project had already been initiated by
Crew \cite{bib:crew1}, \cite{bib:crew2}, who obtained (conditioned
on his conjecture as needed) analogues of several
key results, such as the construction of global monodromy and the theory
of determinantal weights. Our work rests crucially on Crew's initiative.

\subsection{The approaches to Weil II}

To explain our approach to Weil II, it will be helpful to compare and
contrast two proofs in \'etale cohomology: Deligne's original argument
\cite{bib:deligne}, and Laumon's modified version \cite{bib:laumon}.
(See \cite{bib:kw} for more on the Fourier transform method, and its
other applications.)
We begin by describing their common elements, which will also be
common to our approach.

We first describe the basic situation. One begins with a variety
over a finite field $\FF_q$
equipped with a lisse 
sheaf (which we will replace by an overconvergent $F$-isocrystal
in the $p$-adic setting); by induction on dimension, it will suffice to
work on a curve. For each point on the curve, the restriction of the
sheaf to that point is a vector space, on which Frobenius acts via
some linear transformation. Fix an embedding of the algebraic closure
of $\QQ_\ell$ into $\CC$; we define the \emph{weight} of an element of
$\QQ_\ell$ as the base $\sqrt{q}$ logarithm of its complex absolute value.
In this language, we assume that the eigenvalues of Frobenius acting on
each point all have a specific weight, and the goal is to prove that
the eigenvalues of Frobenius on the cohomology of the sheaf have
weights of a certain form. For instance, if we start with the trivial
sheaf, the weights of the eigenvalues of Frobenius at points are all zero,
and the weights of the eigenvalues of Frobenius on cohomology should end up all
being integers less than or equal to 1.

As a first approximation, we consider the \emph{determinantal weights}
of a sheaf; for an irreducible sheaf, these are all equal to the
average of the weights (of the eigenvalues of Frobenius) on cohomology.
One needs to show that these actually coincide with the weights on
cohomology; to begin with, one shows that the determinantal weights
behave correctly under tensor products and the like. (That is, the
determinantal weights of the tensor product of two sheaves are 
pairwise sums of one determinantal weight from each sheaf.) That is 
accomplished by relating determinantal weights to ``global monodromy'',
i.e., the automorphism group of a certain fibre functor (the fibre
at a point) of a certain tannakian category (generated by the given
lisse sheaf). With determinantal weights behaving as predicted,
one then uses an argument about Dirichlet series
(Deligne's analogue of the Rankin squaring trick) to show that they
coincide with weights on cohomology.

The arguments up to this point are essentially ``motivic'', using only
very limited functoriality in cohomology. Thus they can be transcribed
readily into rigid cohomology, as was done by Crew \cite{bib:crew1},
\cite{bib:crew2}. Beyond this point, Deligne's and Laumon's arguments
diverge; we describe them both and then explain how our argument
follows Laumon's method.

In Deligne's original proof, one
passes from the original curve $X$ to $X \times X$, fibers the product
in curves, then resolves singularities of the bad fibres. A careful
analysis of vanishing cycles in this picture makes it possible to
establish that the family has ``large monodromy'', from which one
obtain the desired result. Unfortunately, the appropriate rigid
cohomological analogue of vanishing cycles are not yet completely
understood (by this author, anyway; ongoing work of Crew may shed
light on the situation), 
so it is not so clear how to transcribe
this technique into the rigid setting.

In Laumon's proof (based on a construction of Deligne), one reduces
consideration to the affine line, then performs a geometric analogue of
the Fourier transform. That is, the Fourier transform of a sheaf on
the affine line is a new sheaf, whose fibre at any given point is the
cohomology of the original sheaf twisted by a certain line bundle.
(The operation of twisting and then taking cohomology
is analogous to multiplying a function by a character and then integrating.)
Under certain conditions, the Fourier transform of an irreducible
lisse sheaf can be shown to be irreducible and lisse; this constitutes
another ``large monodromy'' statement that yields the desired result.

It is Laumon's method that we have chosen to imitate here in rigid
cohomology; this approach was suggested by Mebkhout \cite{bib:meb1}.
Specifically, Mebkhout notes that in the analogous context in algebraic
de Rham cohomology, the geometric Fourier transform admits a natural
interpretation in terms of $\calD$-modules, for $\calD$ the ring
of differential operators on the affine line. Namely, the Fourier transform
is essentially the pullback under the automorphism of $\calD$ that
(up to sign) switches multiplication by a coordinate $t$ with differentiation
with respect to $t$. We establish an analogous result for
the Fourier transform of an overconvergent $F$-isocrystal on the affine
line; this is actually a very special case of a general Fourier transform
construction studied in detail by
Huyghe \cite{bib:huyg1}, but we work our special case
out explicitly to keep the discussion self-contained.

One potential difficulty that needs to be skirted here concerns a
present inadequacy in the foundations of $p$-adic cohomology.
Our use of the Fourier transform in rigid cohomology, even in a highly
restricted context, will require us to stray out of the category of
overconvergent $F$-isocrystals, which represent only the constant rank
coefficient objects in rigid cohomology, into the larger category of
arithmetic $\calD$-modules. This theory is expected to carry Grothendieck's
six operations, but progress in this direction is still ongoing.
Fortunately, the particular cohomological operations we will need can
be constructed explicitly, using the pushforward constructions
of \cite{bib:me8}.

It is worth pointing out that a third approach to Weil II is available,
as described by Katz in his lectures
at the 2000 Arizona Winter School \cite{bib:katz}. Katz adopts the
Fourier transform as does Laumon, but avoids a full development of
global monodromy and determinantal weights by explicitly constructing
settings where the global monodromy is ``as large as possible''.
The result is a proof of Weil II more in the manner of Deligne's
first proof of the Weil conjectures \cite{bib:del1}. It should be
possible to adopt Katz's method to the $p$-adic context, but we have
not attempted to do so. (We have however incorporated some of Katz's simplifying
ideas in the context of Laumon's argument.)

\subsection{Loose ends}
\label{subsec:looseends}

The work in this paper by no means constitutes a complete $p$-adic
transcription of the whole story of lisse $\ell$-adic sheaves and their 
cohomology. There are in fact a number of loose ends remaining; we mention
some of them here (with the caveat that some of the descriptions may be
a bit opaque to nonexperts).

Of course, the most notable loose end is the lack of an analogue of
constructible $\ell$-adic sheaves, without which we cannot state Deligne's
Weil II theorem in its relative form. 
Berthelot's theory of arithmetic
$\mathcal{D}$-modules is expected to provide such an analogue, namely
the holonomic arithmetic $\mathcal{D}$-modules. However, there are a number
of stubborn questions that remain outstanding about such objects, hindering
the proof that they admit the expected cohomological operations; see
\cite[Section~5.3.6]{bib:ber7} for a rundown of these.
The ideas of Caro \cite{bib:caro} may be of some use in skirting these
technical difficulties, but this remains to be seen.

Another missing piece of the story is equidistribution of the eigenvalues
of Frobenius, which in the $\ell$-adic situation amounts to the Chebotarev
density theorem. As Richard Crew points out, the result one seeks for
overconvergent $F$-isocrystals fails in the larger but simpler category
of convergent $F$-isocrystals, so some additional ideas may be required.

Also missing is a theorem to the effect that an irreducible
(or geometrically irreducible) overconvergent $F$-isocrystal on a variety
remains so under restriction to a suitable curve. Such a result is needed to
establish that under the hypothesis of ``realizability'', a geometrically
irreducible overconvergent $F$-isocrystal is pure of some weight. Namely,
the Rankin-Selberg method only gives this result on a curve (see
Theorem~\ref{thm:realpure}); in the $\ell$-adic context, one then uses
Bertini to show that a geometrically irreducible sheaf remains so upon
restriction to a suitable curve, so the result follows for any variety.
Again, one cannot prove the analogous $p$-adic results by the same argument,
as it would then apply in the convergent category, whereas a typical
irreducible overconvergent $F$-isocrystal becomes reducible in the convergent
category (because it admits a slope filtration).

Finally, one also expects that one can define local epsilon factors
and prove a product formula for them, following Laumon. This may be addressed
by ongoing work of Marmora; an analogue in the setting of complex local
systems appears in work of Beilinson, Bloch and Esnault.

\subsection{Further directions}

It is worth speculating on what one can actually accomplish with
a $p$-adic form of Weil II. Here are some sample speculations; the reader
is of course encouraged to come up with additional ones.

When using $p$-adic cohomology, one can consider both the archimedean
and $p$-adic valuations of eigenvalues in cohomology
(see Theorem~\ref{thm:padic} for the latter).
This may make it possible to study the interplay between the two,
which may be useful in some geometric settings.
For instance,
if the eigenvalues of Frobenius on an overconvergent $F$-isocrystal
are algebraic integers all having complex absolute value $q^i$ and $p$-adic
valuation $i v_p(q)$, then each one is $q^i$ times a root of unity. That
is, these eigenvalues arise geometrically from Tate classes, which are
expected to be represented by algebraic cycles.

It may also be useful to have a $p$-adic form of Weil II available in
circumstances where one can easily construct a relevant overconvergent
$F$-isocrystal but the ``corresponding'' $\ell$-adic sheaf is not so easy
to find. We suspect examples of this flavor exist in the theory of 
exponential sums, but unfortunately we do not have a concrete one
available to exhibit.

Further down the road,
a $p$-adic form of Weil II may also make it possible to include $p$-adic
cohomology in various ``independence of $\ell$'' assertions. A
distant hope in this direction (which will certainly require resolving
the loose ends of the previous section) is that one can imitate 
Lafforgue's work on the Langlands correspondence for function fields
\cite{bib:lafforgue} in $p$-adic cohomology. For instance, producing
a correspondence from automorphic representations to isocrystals would 
resolve the existence of the ``petit camarade cristalline''
conjecturally associated to a lisse $\ell$-adic sheaf
\cite[Conjecture~1.2.10]{bib:deligne}
(compare \cite[Th\'eor\`eme~VII.6]{bib:lafforgue}).

Finally, $p$-adic Weil II may be relevant in studying forms of the 
weight-monodromy conjecture, which asserts roughly that the monodromy
filtration of the cohomology of a variety over a local field is compatible
with the filtration by weights. See \cite{bib:del1} for the classical form
concerning $\ell$-adic cohomology, and \cite{bib:mokrane} for a variant 
concerning $p$-adic cohomology.

\subsection{Structure of the paper}

We conclude this introduction by describing the contents of
the various chapters of the paper. Before proceeding, 
we caution the reader that this paper
is ``semi-expository'', in that a significant fraction of the assertions
here are not really new. However, some of these assertions
(like the Lefschetz trace formula and the Grothendieck-Ogg-Shafarevich
formula) do
not appear in the literature in the precise form we need them, so
we have decided to err on the side of verbosity.

Chapter~\ref{sec:overcon} is essentially a review of prior constructions
and results on rigid cohomology and overconvergent $F$-isocrystals.
We introduce these notions, point out the connection with
Monsky-Washnitzer cohomology on smooth affine schemes,
and give a trace formula for Frobenius.

In Chapter~\ref{sec:push}, we recall the results of
\cite{bib:me8} on higher direct images in rigid cohomology, along 
some simple morphisms of relative dimension 1.
We then work out some more precise results along these lines,
particularly concerning degeneration in families.
That is, we must understand how the cohomology
of a single member of a family is controlled by the cohomology of
the other members of the family.

In Chapter~\ref{sec:fourier}, we introduce, in a limited context,
the geometric Fourier transform in the $p$-adic setting and its
$\calD$-module interpretation.
We also formulate an analogue of the Grothendieck-Ogg-Shafarevich formula,
which constrains the Euler-Poincar\'e characteristic of (the cohomology) of
an overconvergent $F$-isocrystal in terms of local monodromy.
This formula is needed to show that the Fourier transform of certain
overconvergent $F$-isocrystal are again isocrystals.

In Chapter~\ref{sec:arch}, we quickly summarize Crew's transcription
of the ``common'' parts of
the arguments of Deligne and Laumon.
We then carry out the analogue of Laumon's proof of the Weil II
main theorem. We also give an estimate in the same spirit for
the $p$-adic valuations of eigenvalues in cohomology.

\subsection*{Acknowledgments}
Thanks to SFB 478 ``Geometrische Strukturen in der Mathematik'' at
Universit\"at M\"unster for its hospitality, and to
Vladimir Berkovich, Pierre Berthelot,
Christopher Deninger, Mark Kisin, and Arthur Ogus for helpful discussions. The
author was partially supported by a National Science Foundation
postdoctoral fellowship and by NSF grant DMS-0400727.

\section{Overview of rigid cohomology}
\label{sec:overcon}

We start with a review of rigid cohomology to the extent that we need it
in the calculations of this paper. This is to say that we will not review
very much, deferring instead to \cite{bib:me8} for a more detailed discussion.

Let $q$ be a fixed power of the prime $p$.
Let $k$ be a perfect field of characteristic $p$ containing $\FF_q$, 
let $\gotho$ be a finite totally ramified extension of the ring of
Witt vectors $W(k)$, let $\gothm$ be the maximal ideal of $\gotho$,
and let $K$ be the fraction field of $\gotho$.
We will assume throughout that $\gotho$ admits
an automorphism $\sigma_K$ lifting the $q$-th power map,
which we regard as fixed. For instance, if $\gotho = W(k)$, then there
is a unique choice of $\sigma_K$; it coincides with a power of the Witt
vector Frobenius. Also, if $k = \FF_q$, we may of course take $\sigma_K$
to be the identity map whatever $\gotho$ happens to be.

%We will at times need to enlarge $q$
%(i.e., to study a particular closed point on a curve over $\FF_q$); if
%$q'$ is a power of $q$, we will denote by $\gotho'$ and $K'$ the smallest
%unramified
%extensions of $\gotho$ and $K$, respectively, whose residue
%fields contain $\FF_{q'}$.

We will frequently consider modules over various rings
equipped with a linear or semilinear
endomorphism. If $M$ is such a module equipped with $F$,
we write $M(-i)$ to denote $M$ equipped with $q^i F$,
and call this the \emph{$i$-th Tate twist} of $M$.
%\fixme{check the sign convention everywhere}

\subsection{The formalism of rigid cohomology}
\label{subsec:formalism}

We first recall some of the formalism of
rigid cohomology, following \cite{bib:ber0} and \cite[Chapter~4]{bib:me8}, 
and summarize the key results of
\cite{bib:me8} that we will be using; we postpone defining
anything until the next section. For shorthand, we abbreviate
``separated scheme of finite type over (the field) $k$'' to
``variety over $k$''; for our purposes, 
there would be no real harm to bundling the
adjective ``reduced'' as well.

The coefficient objects in rigid cohomology are called
\emph{overconvergent $F$-isocrystals} (with respect to $K$); 
they form a category fibred in symmetric 
tensor categories over the category of $k$-varieties.
In other words:
\begin{itemize}
\item
The fibre over each variety $X$
 admits direct sums, tensor products (which commute), 
duals, internal Homs, and
an identity object $\calO_X$ for tensoring (the ``constant sheaf''). 
\item
To each morphism $f: X \to Y$ of $k$-varieties is associated a pullback functor
$f^*$
that commutes with the aforementioned operations. These pullback
functors compose up to natural isomorphism.
\end{itemize}
This particular category has the following additional properties:
\begin{itemize}
\item
One also has pullback functors associated to automorphisms of $K$.
\item
The category is equipped with a natural isomorphism $F$ (``Frobenius'')
between the
pullback functor associated to $\sigma_K$ and the identity functor.
(Beware: the action of $F$ on the dual of an overconvergent
$F$-isocrystal is the \emph{inverse} transpose of its action on the
original.)
\item
The fibre over $\Spec k'$, for $k'$ a finite extension of $k$,
is equivalent to the category of finite dimensional $K'$-vector spaces,
for $K'$ the unramified extension of $K$ with residue field $k'$,
equipped with a bijective $\sigma_{K'}$-linear transformation $F$.
(Here $\sigma_{K'}$ is the unique extension of $\sigma_K$
to an automorphism of $K'$ lifting the $q$-power Frobenius.)
\item
There are Tate twist functors which pointwise multiply $F$ by 
the appropriate power of $q$.
\end{itemize}

Associated to a $k$-variety $X$ and an overconvergent $F$-isocrystal
$\calE$ over $X$ are its rigid cohomology spaces
 $H^i_{\rig}(X/K,\calE)$ and its rigid cohomology spaces with compact
supports
$H^i_{c,\rig}(X/K, \calE)$. These are vector spaces over $K$,
which coincide if $X$ is proper;
they vanish for $i<0$ or $i > 2 \dim X$
and are finite dimensional in general by
\cite[Theorems~1.1 and~1.2]{bib:me8}.
If $X$ is smooth of pure dimension $n$, by
\cite[Theorem~1.3]{bib:me8} there is a canonical perfect pairing 
(Poincar\'e duality)
\[
H^i_{\rig}(X/K, \calE) \times H^{2n-i}_{c,\rig}(X/K, \calE^\dual)
\to \calO_X(-n).
\]

The cohomology spaces are functorial in the following senses. 
Given overconvergent $F$-isocrystals
$\calE_1, \calE_2$ on $X$ and a morphism $h: \calE_1 \to \calE_2$,
we obtain morphisms 
\[
H^i_{\rig}(X/K, \calE_1) \to H^i_{\rig}(X/K, \calE_2),
\qquad
H^i_{c,\rig}(X/K, \calE_1) \to H^i_{c,\rig}(X/K, \calE_2)
\]
which compose as expected.
Given a morphism $f:X \to Y$ of varieties
and an overconvergent $F$-isocrystal $\calE$ on $Y$, we obtain
morphisms
\[
H^i_{\rig}(Y/K, \calE) \to H^i_{\rig}(X/K, f^* \calE)
\]
which again compose as expected.
If $f: X \to Y$ is finite \'etale, there is a pushforward functor
$f_*$ from overconvergent $F$-isocrystals on $X$ to those on $Y$,
and we have canonical isomorphisms
\[
H^i_{\rig}(X/K, \calE) \cong H^i_{\rig}(Y/K, f_* \calE), \qquad
H^i_{c,\rig}(X/K, \calE) \cong H^i_{c,\rig}(Y/K, f_* \calE).
\]
Also, cohomology is unchanged by passing from $X$ to its reduced
subscheme.

In cohomology with compact supports, we have an excision exact sequence
\begin{equation} \label{eq:excis2}
\cdots \to H^i_{c,\rig}(X \setminus Z/K, \calE) \to
H^i_{c,\rig}(X/K, \calE) \to H^i_{c,\rig}(Z/K, \calE) \to 
H^{i+1}_{c,\rig}(X \setminus Z/K, \calE) \to \cdots,
\end{equation}
where the maps ``at one level'' (i.e., from one $H^i$ to another)
are Frobenius-equivariant.
There is also an excision sequence in ordinary cohomology:
\begin{equation} \label{eq:excis1}
\cdots \to H^i_{Z,\rig}(X/K, \calE) \to H^i_{\rig}(X/K, \calE)
\to H^i_{\rig}(X \setminus Z/K, \calE) \to H^{i+1}_{Z,\rig}(X/K, \calE)
\to \cdots.
\end{equation}
but it includes the relative cohomology term $H^i_{Z,\rig}(X/K,\calE)$,
which we will not discuss here. (It can be related to ordinary
cohomology via a Gysin map, as in \cite{bib:tsu5}.)

For any closed point $x$ of a variety $X$, we can pull back an
overconvergent $F$-isocrystal $\calE$ along the embedding $x \hookrightarrow
X$ to obtain an object we notate $\calE_x$. As noted above,
the data of $\calE_x$ amounts
to a vector space over the unramified extension $K'$ of
$K$ with residue field $\kappa(x)$, equipped with a $\sigma_{K'}$-linear
bijection induced by $F$. 
We call either object the \emph{fibre} of $\calE$ at $x$.

Now suppose $k = \FF_q$ and $\sigma_K$ is the identity morphism.
Then $F^{\deg(x)}$ induces a \emph{linear} transformation $F_x$ on $\calE_x$.
(However, the natural action of $F^{\deg(x)}$ on $\calE_x \otimes_{K'} L$, for
$L$ a finite extension of $K'$, is typically \emph{not} linear.)
We then have a Lefschetz trace formula for Frobenius, given by the
following equality of formal power series:
\begin{equation} \label{eq:fixpt2}
\prod_{x \in X} \det(1 - F_x t^{\deg(x)}, \calE_x)^{-1}
= \prod_i \det(1 - F t, H^i_{c,\rig}(X/K, \calE))^{(-1)^{i+1}}.
\end{equation}
This is a result of \'Etesse and le Stum 
\cite[Th\'eor\`eme~6.3]{bib:etesse-lestum}; we will briefly review its
derivation in Chapter~\ref{sec:trace}. Note that  in 
\eqref{eq:fixpt2}, the determinant of $1 - F_x t^{\deg(x)}$ is being taken
over $K'$, but actually has coefficients in $K$. If one prefers to work
exclusively over $K$, one may write \eqref{eq:fixpt2} in the form given
in \cite{bib:etesse-lestum}:
\[
\prod_{x \in X} \det_K(1 - F_x t^{\deg(x)}, \calE_x)^{-1/\deg(x)}
= \prod_i \det(1 - F t, H^i_{c,\rig}(X/K, \calE))^{(-1)^{i+1}}.
\]

\subsection{Interlude: rigid cohomology after Berthelot}

In this section, we recall Berthelot's definition of overconvergent
$F$-isocrystals and rigid cohomology in the case of a quasiprojective
variety $X$. (The general case can be obtained by glueing.)
This section is not required reading for the rest of the paper, as
we will use a less general but more convenient construction of
Monsky-Washnitzer for computations (for which skip to the next section);
it is here merely to illustrate the explicit nature of Berthelot's
construction even in fairly wide generality. See 
\cite[Chapter~4]{bib:me8} or \cite{bib:ber0} for more of an overview,
or \cite{bib:ber2} for a discussion with more proofs (at least for
constant coefficients).

Suppose for simplicity that
$X$ is a quasiprojective variety over $k$.
Choose an open immersion $X \hookrightarrow Y$ with $Y$ projective
and a
closed immersion $Y \hookrightarrow \PP^n_k$; put $Z = Y \setminus X$.
Let $P$ be the analytic projective space over $K$; this is a rigid analytic
space equipped with a specialization map $\speci: P \to \PP^n_k$. (By a
construction of Raynaud, one can identify the points of $P$ with 
those closed formal subschemes of the formal completion of
$\PP^n_{\gotho}$ along $\PP^n_k$
which are integral and finite flat over $\gotho$.) Define
the \emph{tubes} $]X[$, $]Y[$, $]Z[$ as the inverse images under
$\speci$ of $X$, $Y$, $Z$,
respectively.

A \emph{strict neighborhood} of $]X[$ in $]Y[$ is an admissible (in the
sense of rigid analytic geometry) 
open subset $V$ of $]Y[$ containing $]X[$, with
the property that for every affinoid $U$ contained in $]Y[$ and any
functions $f_1, \dots, f_d \in \Gamma(U, \calO_U)$ such that
\[
U \cap ]Z[ = \{x \in U: |f_i(x)| < 1 \quad (1 \leq i \leq d) \},
\]
there exists $\lambda < 1$ such that
\[
U \setminus V \subseteq
\{ x \in U: |f_i(x)| < \lambda \quad (1 \leq i \leq d) \}.
\]
Let $\calO^\dagger_{]X[}$ be the direct limit of $j^\dagger \calO_V$ over all
inclusions $j\!\!: ]X[ \hookrightarrow V$ of $]X[$ into a strict neighborhood.
For example, if $X = \AAA^1$ and $Y = \PP^1$, then $]X[$ is the closed
unit disc, whereas $\Gamma(]\overline{X}[, \calO^\dagger_{]X[})$ 
consists of series convergent
on some strictly larger closed disc. Let $\Omega^{i,\dagger}_{]X[}$ be
the direct limit of $j^\dagger \Omega^i_V$ over all $j$ as above.

Choose an extension of $\sigma_K$ to a map from $]X[$ to itself
which reduces to the $q$-power Frobenius under $\speci$.
(This can always be done because $P$ admits such a map.)
Then the category of overconvergent $F$-isocrystals on $X$ is equivalent to the
category of finite locally free $\calO^\dagger_{]X[}$-modules $\calE$,
equipped with an integrable $K$-linear
 connection $\nabla\!: \calE \to \calE \otimes 
\Omega^{1,\dagger}_{]X[}$
and an isomorphism $F: \sigma^* \calE \stackrel{\sim}{\to} \calE$
of modules with connection.
In particular, the latter category does not depend on any of the choices
made so far. (In practice, one \emph{defines} the category of overconvergent
$F$-isocrystals this way, then verifies the independence from
choices.) Note that all of the objects defined with daggers are actually
defined over some strict neighborhood.

The rigid cohomology of $X$ with coefficients in
an overconvergent $F$-isocrystal $\calE$ can be calculated
as the cohomology of the de~Rham complex
\[
H^i_{\rig}(X/K, \calE) = H^i(]X[, \calE \otimes \Omega^{., \dagger}_{]X[}).
\]
We must work a bit harder to compute rigid cohomology with compact supports.
For $V$ a strict neighborhood of $]X[$ on which $\calE$ is defined,
let $\iota\!:\, (]Z[ \cap V) \hookrightarrow V$ be the canonical inclusion.
The cohomology with compact supports will be defined in terms of the
left exact functor on sheaves on abelian groups
\[
\underline{\Gamma}_V(E) = \ker(E \mapsto \iota_* \iota^* E);
\]
for $E$ coherent over $\calO_{V}$, we have $R^q \iota_* \iota^* E = 0$
for $q \geq 1$, so the derived total complex $\mathbb{R} 
\underline{\Gamma}_V(E)$
is isomorphic in the derived category to the two-term complex
$0 \to E \to \iota_* \iota^* E \to 0$.
The cohomology of $X$ with compact
supports with coefficients in $\calE$ is then given by
\[
H^i_{c,\rig}(X/K, \calE) = H^i(V, \mathbb{R} \underline{\Gamma}_{V}(\calE
\otimes \Omega^{.}_{V}));
\]
this turns out to be independent of the choice of $V$.

In both cases, the choice of $K$ is somewhat auxiliary: the computation
of cohomology commutes with replacing $K$ by a finite extension.

\subsection{Affinoid and dagger algebras}
\label{subsec:dagger}

We compute in rigid cohomology using a simplified
mechanism due to Monsky and Washnitzer; this theory looks like
algebraic de Rham cohomology except that the coordinate
ring of the original affine scheme is replaced by a ``dagger algebra''.
In this section, we recall the construction and properties of dagger
algebras, following \cite[Chapter~2]{bib:me8}.

We first recall the notion of an affinoid algebra. Define the ring
\[
T_{n} = K \langle x_1, \dots, x_n \rangle = \left\{ \sum_I a_I x^I: a_I \in K,
\lim_{\sum I \to \infty} |a_I| = 0 \right\}.
\]
Here $I = (i_1, \dots, i_n)$ denotes an $n$-tuple of nonnegative integers,
$x^I = x_1^{i_1}\cdots x_n^{i_n}$, and $\sum I = i_1 + \cdots + i_n$.
An \emph{affinoid algebra} over $K$ is any $K$-algebra isomorphic to
a quotient of $T_n$ for some $n$.
If $A$ is a reduced affinoid algebra,
there is a canonical power-multiplicative norm $|\cdot|_{\sup,A}$ on $A$, 
called the \emph{spectral norm}, with respect to which $A$ is complete.
We also define the \emph{spectral valuation} $v_A$ by 
\[
v_A(x) = -\log_p |x|_{\sup,A}.
\]

We now proceed to dagger algebras.
For $\rho>1$ in the norm group of $K^{\alg}$, define the ring
\[
T_{n,\rho} = \left\{ \sum_I a_I x^I: a_I \in K,
\lim_{\sum I \to \infty} |a_I| \rho^{\sum I} = 0 \right\};
\]
it is an affinoid algebra with spectral norm given by
$|\sum a_I x^I|_{\rho} = \max_I \{|a_I| \rho^{\sum I}\}$.
Define the ring of \emph{overconvergent power series} in $x_1,\dots,x_n$ by
\[
W_n = K\langle x_1, \dots, x_n \rangle^\dagger
= \bigcup_{\rho>1} T_{n,\rho}. 
\]
We note in passing that any finite projective module over $T_n$ or $W_n$ is
free,
by an analogue of the Quillen-Suslin theorem; see
 \cite[Theorem~6.7]{bib:mefull}.
A \emph{dagger algebra} over $K$ is any $K$-algebra isomorphic
to a quotient of $W_n$ for some $n$. Topologizing $W_n$ as a subspace
of $T_n$, we induce a topology on any dagger algebra, called the
\emph{affinoid topology}.

If $A$ is a dagger algebra, we define a \emph{fringe algebra} of $A$
as a subalgebra of the form $f(T_{n,\rho})$ for some surjection
$f: W_n \to A$ and some $\rho>1$ in the norm group of $K^{\alg}$;
note that any fringe algebra is an affinoid algebra, and so has
a natural topology under which it is complete.
We can retopologize $A$ as the direct
limit of its fringe algebras (i.e., a sequence converges to a limit if and
only if it does so in some fringe algebra); we call this topology the
\emph{fringe topology}. The fringe topology is crucial for 
constructing Robba rings over dagger algebras in Section~\ref{subsec:robba2},
and for
obtaining
the Lefschetz trace formula in Chapter~\ref{sec:trace}.

Let $T_n^{\inte}$ or $W_n^{\inte}$ be the subring of $T_n$ or 
$W_n$, respectively, consisting of series with integral
coefficients. Then it turns out that the image of $T_n^{\inte}$ or 
$W_n^{\inte}$
under a surjection $f: T_n \to A$ or $f: W_n \to A$
is independent of $f$; we call it the
\emph{integral subring} of $A$, denoted $A^{\inte}$. More generally,
any homomorphism $g: A \to B$ of affinoid or 
dagger algebras carries $A^{\inte}$
into $B^{\inte}$.

In the same vein, it turns out that the image under a surjection
$f: W_n \to A$ of the ideal of $W_n^{\inte}$ consisting
of series whose coefficients all lie in $\gothm$
is independent of $f$. The elements of this ideal are the
topologically nilpotent elements of $A^{\inte}$;
the quotient of $A^{\inte}$ by this ideal, which is finitely generated
as a $k$-algebra, is called the \emph{reduction} of $A$.
If $R$ is the reduction of $A$, we call $\Spec R$ the 
\emph{special fibre} of $A$.

Given a dagger algebra $A = W_n/\gotha$, write
\[
A \langle t \rangle^\dagger = W_{n+1}/\gotha W_{n+1},
\]
identifying $t$ with $x_{n+1}$. This construction does not depend on
the presentation of $A$. For $f \in A$, write
\[
A \langle f^{-1} \rangle^\dagger = A \langle t \rangle^\dagger / (tf - 1);
\]
this is called the \emph{localization} of $A$ at $f$.

\subsection{Cohomology of affine schemes}
\label{subsec:mw}

We now construct Monsky-Washnitzer cohomology, our main computational
tool in studying rigid cohomology on smooth affine varieties.
Our reference now is \cite[Chapter~3]{bib:me8}.

The \emph{module of continuous differentials} $\Omega^1_{A}$ of a dagger
algebra can be constructed as follows. For $A = W_n$, take it to be the
free module generated by $dx_1, \dots, dx_n$ equipped with the $K$-linear
derivation $d: W_n \to \Omega^1_{W_n}$ given by
\[
\sum_I c_I x^I \mapsto \sum_I \sum_{j=1}^n i_j c_I (x^I/x_j)\,dx_j.
\]
For $A \cong W_n/\gotha$, let $\Omega^1_{A}$ be the quotient of
$\Omega^1_{W_n} \otimes_{W_n} A$ by the submodule generated by
$dr$ for $r \in \gotha$. This construction ends up being universal for
$A$-linear derivations into finitely generated $A$-modules; in particular,
it yields a well-defined $A$-module
$\Omega^1_{A}$ and $K$-linear derivation $d: A \to \Omega^1_{A}$.
If $A$ is a subring of the dagger algebra $B$, we define the relative
module of differentials $\Omega^1_{B/A}$ as the quotient of
$\Omega^1_{B}$ by the images of $da$ for $a \in A$.
We also put $\Omega^i_{B/A} = \wedge^i_A \Omega^1_{B/A}$.

At this point, we restrict to a special class of dagger algebras. We say
a dagger algebra $A$ is \emph{of MW-type} if the ideal of topologically
nilpotent elements of $A^{\inte}$ is generated by a uniformizer of $\gotho$
and the special fibre of $A$ is smooth.
In the terminology of \cite{bib:mw1},
$B$ is a formally smooth, weakly complete, weakly finitely generated algebra
over $(\gotho, \gothm)$.

A \emph{Frobenius lift} on a dagger algebra $A$ of MW-type is a 
ring endomorphism $\sigma: A \to A$ acting on $K$ via $\sigma_K$
and acting on $A^{\inte}
\otimes_\gotho k$ as the $q$-th power map $x \mapsto x^q$.
Such a map exists for any $A$;
for example, if $A = W_n$, we can define a \emph{standard Frobenius}
$\sigma$ by the formula
\[
\left( \sum_I c_I t^I \right)^\sigma = \sum_I c_I^{\sigma_K} t^{qI}.
\]

Given a dagger algebra $A$ equipped with a Frobenius lift $\sigma$,
we define a \emph{$\sigma$-module} over $A$ as a finite locally free
$A$-module equipped with
\begin{enumerate}
\item[(a)] a \emph{Frobenius structure}:
an additive, $\sigma$-linear map $F: M \to M$ (that is, $F(a \bv)
= a^\sigma F(\bv)$ for $a \in A$ and $\bv \in M$) which induces
an isomorphism $\sigma^* M \to M$.
\end{enumerate}
We define a \emph{$(\sigma, \nabla)$-module} over $A$
as a $\sigma$-module additionally equipped with
\begin{enumerate}
\item[(b)] an \emph{integrable connection}:
an additive, $K$-linear map $\nabla: M \to M \otimes_A
\Omega^1_{A}$ satisfying the Leibniz rule: $\nabla(a \bv) = a
\nabla(\bv) + \bv \otimes da$ for $a \in A$ and $\bv \in M$, and
such that, if we write $\nabla_n$ for the induced map $M \otimes_A
\Omega^n_{A} \to M \otimes_A \Omega^{n+1}_{A}$, we have
$\nabla_{n+1} \circ \nabla_n = 0$ for all $n \geq 0$;
\end{enumerate}
subject to the compatibility condition
\begin{enumerate}
\item[(c)] the isomorphism $\sigma^* M \to M$ induced by $F$ is horizontal
for the corresponding connections; in other words, the following diagram 
commutes:
\[
\xymatrix{
M \ar^{F}[d] \ar^(.3){\nabla}[r] & M \otimes_A \Omega^1_{A}
\ar^{F \otimes d\sigma}[d] \\
M \ar^(.3){\nabla}[r] & M \otimes_A \Omega^1_{A}.
}
\]
\end{enumerate}

For example, the module $M = A$, with $F$ acting by $\sigma$ and
$\nabla$ acting by $d$, is a $(\sigma, \nabla)$-module, called
the \emph{trivial $(\sigma, \nabla)$-module}. More generally, if
$M$ is spanned over $A$ by the kernel of $\nabla$, we say
$M$ is \emph{constant}.

Given a $(\sigma, \nabla)$-module $M$ over $A$, we define the cohomology
spaces as the cohomology of the de Rham complex tensored with $M$. That is,
\[
H^i(M) = \frac{\ker(\nabla_{i}: M \otimes_A \Omega^i_{A} \to
M \otimes_{A} \Omega^{i+1}_{A})}{\im(\nabla_{i-1}: M \otimes_{A}
\Omega^{i-1}_{A}
\to M \otimes_A \Omega^i_{A})}.
\]

If $M$ is a $(\sigma, \nabla)$-module over $A$, 
we call an $A$-submodule $N$ of $M$
a \emph{$(\sigma, \nabla)$-submodule} if it is closed under
$F$ and $\nabla$ (the latter meaning that $\nabla(N) \subseteq
N \otimes \Omega^1_A$); it turns out \cite[Lemma~3.3.4]{bib:me8}
 that this forces
$N$ to be a direct summand of $M$ as an $A$-module, so the quotient
$M/N$ is also a $(\sigma, \nabla)$-module. (Beware that $N$ need not be
a direct summand of $M$ in the category of $(\sigma, \nabla)$-modules over
$A$; that is, the exact sequence $0 \to N \to M \to M/N \to 0$ may not
have a horizontal splitting.)
This gives us a notion of irreducibility of a $(\sigma, \nabla)$-module; we
say $M$ is \emph{absolutely irreducible} if it remains irreducible
whenever we replace $q$ by a power $q'$ of $q$,
$k$ by a finite extension $k'$ containing $\FF_{q'}$, and $K$ by a
finite extension $K'$ with residue field $k'$.

%We can extend the notion of $(\sigma, \nabla)$-module, and the
%associated definitions in this section, to any pair
%of rings $R$ and $A$ given an endomorphism $\sigma$ on $R$ preserving $A$
%and a definition of $\Omega^1_{R/A}$ such that $d\sigma: \Omega^1_{R/A}
%\to \Omega^1_{R/A}$ is well-defined. In particular, once we define
%the Robba ring $\calR_A$ over $A$, we will have the notion of 
%a $(\sigma, \nabla)$-module over $R = \calR_A$ relative to $A$.

We now summarize 
the relationship of this construction to rigid cohomology;
see \cite[Proposition~1.10]{bib:ber2} for more details
in the constant coefficient case (the general case is similar).
If $X$ is a smooth affine $k$-variety and $\calE$ is an overconvergent
$F$-isocrystal on $X$, then $\calE$ can be identified with a module $M$
over the ring of functions on $]X[$ which extend to some strict 
neighborhood, equipped with a connection $\nabla$. That ring of
functions is a dagger algebra $A$ with special fibre $X$, and given
any Frobenius lift $\sigma$ on $A$, the isomorphism $F: \sigma^* \calE
\stackrel{\sim}{\to} \calE$ gives a Frobenius structure on $M$. In 
other words, $M$ is a $(\sigma, \nabla)$-module over $A$,
and there is a canonical isomorphism
\[
H^i_{\rig}(X/K, \calE) \cong H^i(M).
\]
In particular, the category of $(\sigma, \nabla)$-modules over $A$
is canonically independent of the choice of $\sigma$.

\subsection{The Robba ring and $p$-adic local monodromy}
\label{subsec:robba1}

We next want to make more explicit the cohomology of curves, but first
we need to introduce an auxiliary ring from the theory of $p$-adic
differential equations.

The \emph{Robba ring} $\calR_K = \calR^t_K$ (the latter notation being used
when we need to name the series parameter)
is defined as the ring of bidirectional
power series $\sum_{n=-\infty}^\infty c_n t^n$, with $c_n \in K$, such that
for $r>0$ sufficiently small (depending on the series),
\[
\lim_{n \to \pm \infty} (v_p(c_n) + rn) = \infty.
\]
That is, such a series converges for $t \in K^{\alg}$ satisfying
$\eta < |t| < 1$, for some $\eta$ depending on the series.

We denote by $\calR_K^{\inte}$ the subring of $\calR_K$ of series with
$v_p(c_n) \geq 0$ for all $n$, and by $\calR^+_K$ the subring of series
with $c_n = 0$ for $n<0$. We denote by $\calR_K^{+,\inte}$ the intersection
of these two subrings; it coincides with $\gotho\llbracket t \rrbracket$.

Given $r>0$ rational, for those elements $x = \sum c_n t^n
\in \calR_K$ for which $v_p(c_n) + rn \to \infty$ as $n \to \pm \infty$,
we put
\[
w_r(x) = \inf_n \{v_p(c_n) + rn\};
\]
this function is a discrete valuation on the subring where it is defined. 
Note that for any fixed $x \in \calR_K$, $w_r(x)$ 
is defined for all sufficiently small $r>0$.

We define $(\sigma, \nabla)$-modules over $\calR_K$ or $\calR_K^+$ as
in the dagger algebra setting, replacing $\Omega^1$ with the free module
generated by $dt$. Note that finite locally free modules 
over $\calR_K$ or $\calR_K^+$ are automatically free, because a theorem of 
Lazard \cite{bib:lazard} implies that  $\calR_K$ and $\calR^+_K$ are 
 B\'ezout rings (rings in which every finitely
generated ideal is principal).

A technique due to Dwork (analytic continuation via Frobenius)
leads to the following result; see \cite[Lemma~6.3]{bib:dej1} for its proof.
\begin{lemma} \label{lem:dworktrick}
Let $M$ be a $(\sigma, \nabla)$-module over $\calR_{K}^+$. Then
there exists a basis $\bw_1, \dots, \bw_n$ of $M$ such that
$\nabla \bw_i = 0$ for each $i$. (Note that on any such basis,
$F$ acts via a matrix over $K$.)
\end{lemma}

A weaker form of Lemma~\ref{lem:dworktrick} holds for $M$ over $\calR_K$,
but is much deeper. It is the so-called ``$p$-adic local monodromy theorem'',
and underpins this entire article as well as \cite{bib:me8}.
Proofs have been given by Andr\'e \cite{bib:andre},
Mebkhout \cite{bib:meb2}, and the author \cite{bib:me7}.
\begin{prop} \label{prop:locmono}
Let $M$ be a $(\sigma, \nabla)$-module over $\calR_K$. Then there
exist a finite \'etale extension $\calR'$ of $\calR^{\inte}_K$
and a basis $\bw_1, \dots, \bw_n$ of 
$M \otimes_{\calR^{\inte}_K} \calR'$ such that for $i=1,\dots, n$,
the span $M_i$ of $\bw_1, \dots, \bw_i$ is carried into $M_i \otimes dt$
by $\nabla$
and the image of $\bw_i$ in $M_i/M_{i-1}$ is killed by $\nabla$.
\end{prop}

We say $M$ is \emph{unipotent} if it satisfies the conclusion of the
$p$-adic local monodromy theorem with $\calR' = \calR^{\inte}_K$.
In that case, one can find a basis $\bv_1, \dots, \bv_n$ of $M$ whose
$K$-span is preserved by the operator $E: M \to M$ defined by
$\nabla(\bv) = E(\bv) \otimes \frac{dt}{t}$ 
(see \cite[Proposition~5.2.6]{bib:me8}).

%It is easy to give a multidimensional generalization of the Robba
%ring; this generalization will be useful in Chapter~\ref{sec:push}.
%For $T$ a finite set,
%let $\calR^T_K$ be the set of formal sums $\sum_I c_I x^I$
%with $I$ running over $\ZZ^T$ (i.e., tuples of integers indexed by the
%elements of $T$), such that
%for some $r>0$ depending on the series,
%\[
%\lim_{I \to \infty} v_p(c_I) + \sum_{t \in T} s_t i_t = \infty
%\]
%for $0 < s_t \leq r$. Here the limit means that the left side is below
%any particular cutoff for only finitely many $I$.

\subsection{Cohomology of curves}
\label{subsec:cohcurve}

A detailed study of the cohomology of overconvergent $F$-isocrystals
on curves has been made by Crew \cite{bib:crew2}; we summarize his
results in this section. (Note that Crew's hypothesis of ``strictness''
is superfluous in this setting, thanks to the $p$-adic local monodromy
theorem.)

Let $X$ be a smooth irreducible affine curve, let $\overline{X}$ be its
smooth compactification, and let $A$ be a dagger algebra of MW-type with
special fibre $X$.
Then for each closed point $x \in \overline{X}$, one gets
a (noncanonical)
embedding $A \hookrightarrow \calR_x$,
where $\calR_x$ is a copy of the
Robba ring over the unramified extension $K'$ of $K$ with
residue field $\kappa(x)$; we can and will take this embedding to map into
$\calR_x^{+}$ if $x \in X$. Observe that given such an embedding,
any Frobenius lift on $A$ can be extended
compatibly to $\calR_x$.

%Let $\calE$ be an overconvergent $F$-isocrystal on $X$. Choose
% a Frobenius lift $\sigma$ on $A$, and let $M$ be the
%$(\sigma, \nabla)$-module
%over $A$ corresponding to $\calE$.
%For $x \in X$, apply Lemma~\ref{lem:dworktrick} to $M$ over
%$\calR_x^+$ and let
%$M_x$ be the $K'$-vector space spanned by the resulting basis
%$\bw_1, \dots, \bw_n$.
%Then $M_x$, as a vector space equipped with a $\sigma$-linear
%endomorphism, is isomorphic to the fibre $\calE_x$.

Define
\begin{gather*}
A_{\loc} = \bigoplus_{x \in \overline{X} \setminus X}
\calR_x, \qquad \Omega^1_{\loc} = \Omega^1_A \otimes_A A_{\loc} \\
A_{\qu} = A_{\loc}/A, \qquad \Omega^1_{\qu} = \Omega^1_A \otimes_A A_{\qu}
= \Omega^1_{\loc}/\Omega^1_A,
\end{gather*}
where the last equality holds because $\Omega^1_A$ is a flat $A$-module.
(Note that $A_{\loc}$ is a ring but $A_{\qu}$ is only an $A$-module.)
For $M$ a $(\sigma, \nabla)$-module over $A$
corresponding to an overconvergent $F$-isocrystal $\calE$
on $X$, we have already defined
\begin{align*}
  H^0(M) &= \ker(\nabla: M \to M \otimes_A \Omega^1_A) \\
H^1(M) &= \coker(\nabla: M \to M \otimes_A \Omega^1_A),
\end{align*}
and observed that $H^i(M) \cong H^i_{\rig}(X/K, \calE)$.
We now define
\begin{align*}
  H^0_{\loc}(M) &= \ker(\nabla: M \otimes_A A_{\loc} \to M \otimes_A \Omega^1_{\loc}) \\
  H^1_{\loc}(M) &= \coker(\nabla: M \otimes_A A_{\loc} \to M \otimes_A \Omega^1_{\loc}) \\
H^1_c(M) &= \ker(\nabla: M \otimes_A A_{\qu} \to M \otimes_A \Omega^1_{\qu}) \\
H^2_c(M) &= \coker(\nabla: M \otimes_A A_{\qu} \to M \otimes_A \Omega^1_{\qu});
\end{align*}
Crew \cite{bib:crew2} has shown
that $H^i_c(M) \cong H^i_{c, \rig}(X/K, \calE)$.
(This identification and the previous one become $F$-equivariant once we
specify that $F$ acts on $\Omega^1$ via the linearization $d\sigma$
of the Frobenius lift.)
For $x \in \overline{X} \setminus X$, we write
$H^0_{\loc,x}(M)$ for the kernel of $\nabla: M \otimes_A \calR_x
\to M \otimes_A \Omega^1_{\calR_x}$, so that
$H^0_{\loc}(M) = \oplus_x H^0_{\loc,x}(M)$;
we also write $H^i_{\loc}(X/K, \calE)$ for $H^i_{\loc}(M)$.

All of the $H^i(M)$, $H^i_{\loc}(M)$, and $H^i_c(M)$
are finite dimensional vector spaces over $K$,
by
\cite[Theorem~9.5 and Proposition~10.2]{bib:crew2} and the $p$-adic
local monodromy theorem.
Because the rows of the diagram
\[
\xymatrix{
0 \ar[r] & M \ar[r] \ar[d] & M \otimes_A A_{\loc} \ar[r] \ar[d] & M \otimes_A A_{\qu} \ar[r] \ar[d] & 0 \\
0 \ar[r] & M \otimes_A \Omega^1_{A} \ar[r] & M \otimes_A 
\Omega^1_{\loc} \ar[r] & M \otimes_A \Omega^1_{\qu} \ar[r] & 0
}
\]
are exact, the snake lemma produces the canonical exact sequence
\begin{equation} \label{eq:snake}
0 \to H^0(M) \to H^0_{\loc}(M) \to H^1_c(M)
\to H^1(M) \to H^1_{\loc}(M) \to H^2_c(M) \to 0.
\end{equation}
Moreover, there are $F$-equivariant perfect pairings
\[
H^i(M) \times H^{2-i}_c(M^\dual) \to H^2_c(K) = K(-1)
\]
which correspond to Poincar\'e duality of overconvergent $F$-isocrystals;
there is also an $F$-equivariant perfect pairing
\[
H^0_{\loc}(M) \times H^{1}_{\loc}(M^\dual) \to K(-1).
\]
We will consider these further in Section~\ref{subsec:pushsupp}.

\section{Pushforwards in rigid cohomology}
\label{sec:push}

The notion of a pushforward (or more precisely, of
higher direct images) is the relative version of the notion of
the cohomology of a single space. 
Picking up a thread from \cite[Chapter~7]{bib:me8}, we consider some
simple pushforwards in relative dimension 1.

%Before considering relative dimension 1, we record here a simple but
%useful fact in a case of relative dimension 0.
%\begin{prop} \label{prop:etalepush}
%If $f: X \to Y$ is a finite \'etale morphism of smooth affine varieties and
%$\calE$ is an overconvergent $F$-isocrystal on $X$, then there is a
%canonical isomorphism $H^i(X, \calE) \cong H^i(Y, f_* \calE)$ for all $i$.
%\end{prop}
%In fact, the underlying de~Rham complexes are isomorphic: if
%$R/S$ is an unramified extension of dagger algebras, we have
%\[
%R \otimes_S \Omega^1_{S/K} \cong \Omega^1_{R/K}.
%\]

\subsection{Robba rings over dagger algebras}
\label{subsec:robba2}

For the calculations in this chapter, we need to extend the definition
of the Robba ring by allowing coefficients not just in $K$, but in a more
general dagger algebra. The correct procedure for doing this is given in
\cite[Section~2.5]{bib:me8}; we quickly review it here.

For $A$ a reduced dagger algebra,
the \emph{Robba ring} $\calR_A = \calR^t_A$ 
is defined as the ring of bidirectional
power series $\sum_{n=-\infty}^\infty c_n t^n$, with $c_n \in A$, such that
for $r>0$ sufficiently small (depending on the series),
$c_n p^{\lfloor rn \rfloor} \to 0$ as $n \to \pm \infty$
in the fringe topology of $A$
(that is, within some fringe algebra depending on $r$).
By \cite[Corollary~2.5.5]{bib:me8}, it is equivalent to require that
\[
\lim_{n \to \pm \infty} v_A(c_n) + rn = \infty
\]
for $r>0$ sufficiently small and that
$c_n p^{\lfloor rn \rfloor} \to 0$
in the fringe topology of $A$ for \emph{one} value of $r$.

%For example, if $A = K \langle x \rangle^\dagger$, the Robba ring
%$\calR^t_A$ consists of series in $t,t^{-1},x$ which, for some
%$\delta>1$ and $r>0$, converge for 
%\[
%|x| \leq \delta, \qquad
%\delta^{-1} \leq t < \min\{|x|^{r}, 1\}.
%\]

We define $\Omega^1_{\calR_A/A}$ as the free module over $\calR_A$
generated by $dt$, equipped with the derivation
\[
d: \calR_A \to \Omega^1_{\calR_A/A} \qquad
\sum_i c_i t^i \mapsto \sum_i i c_i t^{i-1}\,dt.
\]
A quick calculation \cite[Proposition~3.1.4]{bib:me8} shows that the
kernel and cokernel of this derivation are isomorphic to $A$ in the
expected manner. In particular, we 
define the residue map $\Res: \Omega^1_{\calR_A/A} \to A$ by
sending $\sum_i c_i t^i\,dt$ to $c_{-1}$; then 
$\omega \in \Omega^1_{\calR_A/A}$
is in the image of $d$ if and only if $\Res(\omega) = 0$.

We define $(\sigma, \nabla)$-modules over $\calR_A$ (or $\calR^+_A$)
relative to $A$
as expected, using the relative module of differentials $\Omega^1_{\calR_A/A}$
and requiring that the connection $\nabla$ be $A$-linear. Then
the Dwork trick admits the following relative version.
\begin{prop} \label{prop:reldwork}
Let $M$ be a free $(\sigma, \nabla)$-module over $\calR^+_A$ relative
to $A$. Then there exists a basis $\bv_1, \dots, \bv_n$ of $M$ such that
$\nabla \bv_i = 0$ for $i=1, \dots, n$. (Note that on any such basis, $F$ acts
via a matrix over $A$.)
\end{prop}
\begin{proof}
Choose any basis $\be_1, \dots, \be_n$ of $M$, and define the $n 
\times n$ matrix $N$ over $\calR^+_A$ by
\[
\nabla \be_j = \sum_i N_{ij} \be_i \otimes dt.
\]
Write $N = \sum_{l=0}^\infty N_l t^l$; 
then a straightforward induction shows that there is a unique invertible 
$n \times n$ matrix $U = I + \sum_{l=1}^\infty U_l t^l$ over 
$A \llbracket t \rrbracket$ such that $NU + \frac{dU}{dt} = 0$. 
Namely, for each $l>0$, we have
\begin{equation} \label{eq:dwork}
l U_l + \sum_{i=0}^{l-1} N_i U_{l-1-i} = 0
\end{equation}
and this lets us solve for $U_l$ in terms of the $U_i$ for $i<l$.

We next show that $U$ has entries in $\calR^+_A$.
Thanks to the usual Dwork's trick (Lemma~\ref{lem:dworktrick}),
all we must verify is that $U_{l} p^{\lfloor rl \rfloor} \to 0$
within some fringe algebra of $A$ for some $r>0$.
Since $N$ has entries in $\calR^+_A$, for some $s>0$
we can choose a fringe algebra 
$B$ such that $v_B(N_l) + sl \to \infty$ as $l \to \infty$.
By taking $s$ large enough, we can ensure that in fact
$v_B(N_l) + s(l+1) > 0$ for all $l$. Then \eqref{eq:dwork} implies 
easily that
$v_B(l! U_l) + sl > 0$ for all $l$, and so
$v_B(U_l) + rl \to \infty$ for $r > s + 1/(p-1)$.
Therefore $U$ indeed has entries in $\calR^+_A$.

The same argument applied to the basis of $M^\dual$ dual to 
$\be_1, \dots, \be_n$ shows that
the inverse transpose of $U$ has entries in $\calR^+_A$. Consequently the
elements $\bv_1, \dots, \bv_n$ of $M$ defined by
\[
\bv_j = \sum_i U_{ij} \be_i
\]
form a basis with the desired property.
\end{proof}

There is also a relative version of the $p$-adic local monodromy theorem
\cite[Theorem~5.1.3]{bib:me8}, which underlies the pushforward construction
of the next section; however, we will not use it explicitly.

\subsection{Pushforwards with and without supports}
\label{subsec:pushsupp}

Let $X$ be a smooth irreducible $k$-variety, let $\calE$ be an overconvergent
$F$-isocrystal on $\AAA^1 \times X$, and let $f: \AAA^1 \times X
\to X$ denote the implicit
 projection. In \cite[Chapter~7]{bib:me8} are constructed
``generic'' higher direct images $R^i f_* \calE$ and $R^i f_! \calE$
of $f$ over an open
dense subscheme of $X$; of course this is the best one can do
within a category of locally free modules, since the rank of the
corresponding cohomology space may jump at particular fibres.
We now review this construction, which follows the setup of
\cite{bib:crew2} as presented earlier in Section~\ref{subsec:cohcurve}.

Let $A$ be a dagger algebra of MW-type with special fibre $X$,
and let $M$ be a $(\sigma, \nabla)$-module over  $A \langle x
\rangle^\dagger$ corresponding to $\calE$. 
Then we get a map
\[
\nabla_v: M \to M \otimes \Omega^1_{A \langle x \rangle^\dagger} 
\to M \otimes \Omega^1_{A \langle x \rangle^\dagger/A}
\]
from the projection $\Omega^1_{A \langle x \rangle^\dagger} \to 
\Omega^1_{A \langle x \rangle^\dagger/A}$.

Embed $A \langle x \rangle^\dagger$ 
into $\calR_A = \calR^t_A$ 
by mapping $\sum c_i x^i$ to $\sum c_i t^{-i}$.
By analogy with the notations of Section~\ref{subsec:cohcurve},
we put
\[
M^{\loc} = M \otimes_{A \langle x \rangle^\dagger} \calR_A, \qquad
M^{\qu} = M^{\loc} / M
\]
and let
\begin{align*}
\nabla_v^{\loc} &: M^{\loc} \to M^{\loc} \otimes_{\calR_A} 
\Omega^1_{\calR_A/A} \\
\nabla_v^{\qu} &: M^{\qu} \to (M^{\loc} \otimes_{\calR_A} \Omega^1_{\calR_A/A})/(M \otimes_{A \langle x \rangle^\dagger}
 \Omega^1_{A \langle x \rangle^\dagger/A})
\end{align*}
be the maps induced by $\nabla_v$.
We then define
\begin{gather*}
R^0 f_* \calE = \ker(\nabla_v), \qquad
R^1 f_* \calE = \coker(\nabla_v) \\
R^0_{\loc} f_* \calE = \ker(\nabla_v^{\loc}), \qquad
R^1_{\loc} f_* \calE = \coker(\nabla_v^{\loc}) \\
R^1 f_! \calE = \ker(\nabla_v^{\qu}), \qquad
R^2 f_! \calE = \coker(\nabla_v^{\qu});
\end{gather*}
we take $R^i f_* \calE$, $R^i_{\loc} f_* \calE$, $R^i f_! \calE$ to be zero
for values of $i$ not covered by the above list. We also sometimes
write $M$ in place of $\calE$ in this notation.

By the snake lemma, we have an $F$-equivariant exact sequence of $A$-modules
\begin{equation} \label{eq:snake2}
0 \to R^0 f_* \calE \to R^0_{\loc} f_* \calE \to R^1 f_! \calE
\to R^1 f_* \calE \to R^1_{\loc} f_* \calE \to R^2 f_! \calE \to 0.
\end{equation}
Moreover, there are canonical $A$-linear, $F$-equivariant 
Poincar\'e duality pairings
\begin{gather}
R^i f_* \calE \times R^{2-i} f_! \calE^\dual \to A(-1) \label{eq:relpair1}
\\
R^i_{\loc} f_* \calE \times R^{1-i}_{\loc} f_* \calE^\dual \to A(-1)
\label{eq:relpair2}
\end{gather}
obtained from the canonical pairing $[\cdot, \cdot]:
\calE \times \calE^\dual \to
A \langle x \rangle^\dagger$ 
and the residue map $\Res: \Omega^1_{\calR_A/A} \to A$. 

By \cite[Theorem~7.3.2 and Proposition~8.6.1]{bib:me8}, 
we have the following result.
\begin{theorem} \label{thm:pushfwd}
There exists a localization $B$ of $A$ such that
$R^i f_* M_B$, $R^i_{\loc} f_* M_B$,
$R^i f_! M_B$ are overconvergent $F$-isocrystals for all $i$
(where $M_B = M \otimes B \langle x \rangle^\dagger$),
and the Poincar\'e duality pairings are perfect.
\end{theorem}

One can relate the cohomology of an overconvergent $F$-isocrystal
to that of its pushforwards (by a Leray spectral sequence); 
the particular instance of this
relationship that we need
is precisely \cite[Proposition~7.4.1]{bib:me8}.
\begin{prop} \label{prop:leray}
Let $X$ be a smooth irreducible affine $k$-variety,
let $f: \AAA^1 \times X \to X$ be the canonical projection, and let
$\calE$ be an overconvergent $F$-isocrystal on $\AAA^1 \times X$ for which
$R^0 f_* \calE$, $R^1 f_* \calE$ are overconvergent $F$-isocrystals.
Then there are canonical, $F$-equivariant exact sequences
\[
H^i_{\rig}(X/K, R^0 f_* \calE) \to 
H^i_{\rig}(\AAA^1 \times X/K, \calE) \to H^{i-1}_{\rig}(X/K, R^1 f_* \calE)
\]
for each $i$.
\end{prop}
These short exact sequences actually come from a long exact sequence, 
but the connecting maps are not $F$-equivariant (they are off by a Tate
twist).

\subsection{Degeneration in families}

Our strategy for studying the cohomology of an isocrystal on a curve is
to embed that isocrystal into a family most of whose fibres are easy
to control. For this to return a result on the original isocrystal, we
need a theorem that specifies how the cohomology of an isocrystal behaves under
specialization. 
A corresponding statement in \cite{bib:katz} is the ``degeneration
lemma''.

We will need to work over a certain auxiliary ring. (This ring must be
chosen carefully; in a preliminary version of this paper, the wrong auxiliary
ring was used. See Remark~\ref{rem:degen}.)
Let $\calS$ denote the
ring of power series $\sum_{i,j \in \ZZ} c_{ij} s^i t^j$ over $K$
in two variables $s$ and $t$ with the following property: for each $0< \delta < 1$
sufficiently close to 1, there exists $0 < \epsilon < 1$ such that the series
converges for $s,t \in K^{\alg}$ with $|s|= \delta$ and $\epsilon < |t| < 1$.
We use the superscripts $s+$ and $t-$ to denote the restriction of $\calS$
to the subrings where $s$ occurs only with positive powers and where
$t$ occurs only with negative powers, respectively.

The value of the ring $\calS$ is that it is defined using a very mild
convergence restriction on series, so many other rings naturally embed into it.
Specifically, we can and will identify $\calR^{s,+}_{K \langle x \rangle^\dagger}$
and $\calR^s_{K \langle x \rangle^\dagger}$ with subrings of
$\calS^{s+,t-}$ and $\calS^{t-}$, respectively, by identifying $x$ with $t^{-1}$.
In particular, this gives an embedding of $A \langle x \rangle^\dagger$ into
$\calS$ for any localization $A$ of $K \langle s \rangle^\dagger$,
since $A \langle x \rangle^\dagger 
\subset \calR^s_{K \langle x \rangle^\dagger}$.

Given a $(\sigma, \nabla)$-module $M$ (necessarily free)
over $K \langle s,x \rangle^\dagger$,
for any localization $A$ of $K \langle s \rangle^\dagger$, we write
$M_A$ for $M \otimes A \langle x \rangle^\dagger$. We also write
$\nabla_s$ and $\nabla_x = \nabla_t$ 
for the components of $\nabla$ mapping to
$M \otimes ds$ and $M \otimes dx$, respectively.

\begin{theorem} \label{thm:degen}
Let $M$ be a $(\sigma, \nabla)$-module over $K\langle s,x
\rangle^\dagger$, and let $f: K \langle s \rangle^\dagger \to
K \langle s,x \rangle^\dagger$ denote the canonical inclusion.
Let $A$ be a localization of $K \langle s
\rangle^\dagger$ such that the conclusion of Theorem~\ref{thm:pushfwd}
holds for $M_A$ and $M^\dual_A$.
Then there is a canonical $F$-equivariant injection
\[
H^1_c(M/sM) \hookrightarrow H^0_{\loc,s=0}(R^1 f_! M_A).
\]
\end{theorem}
\begin{proof}
By Poincar\'e duality, it is equivalent to exhibit a canonical $F$-equivariant
pairing
\begin{equation} \label{eq:canfeq}
H^1_c(M/sM) \times H^1_{\loc,s=0}(R^1 f_* M_A^\dual) \to K(-2)
\end{equation}
which is nondegenerate on the left.
Using the relative Dwork's trick (Proposition~\ref{prop:reldwork}), 
we get a $K$-linear
map $g: M/sM \to M \otimes \calR^{s,+}_{K \langle x \rangle^\dagger}$
such that for all $\bv \in M/sM$, 
$g(\bv)$ reduces to $\bv$ modulo $s$,
$\nabla_s g(\bv) = 0$, and
$\nabla_t g(\bv) = g(\nabla_t \bv)$.
In particular, $g$ induces an $F$-equivariant inclusion 
\begin{equation} \label{eq:feqincl}
H^1_c(M/sM) \hookrightarrow \frac{ \{ \bv \in M \otimes 
\calS: \quad
\nabla_s \bv = 0,  \quad \nabla_t \bv \in M \otimes \calS^{t-} \otimes dx
\}} {\{\bv \in M \otimes \calS^{t-}: \quad \nabla_s \bv = 0\}}.
\end{equation}

We may as well assume that $s^{-1} \in A$. Then
we have a natural $F$-equivariant perfect pairing
\[
\frac{M \otimes \calS}{M \otimes \calS^{t-}} \times (M_A^\dual \otimes 
(ds \wedge dx)) \to K(-2)
\]
given by the residue map on $\calS \otimes (ds \wedge dt)$ (i.e., extracting
the coefficient of $(ds/s) \wedge (dt/t)$).
If we restrict on the left to the classes of those
$\bv$ with $\nabla_t \bv \in M \otimes \calS^{t-} \otimes dx$, 
then the pairing vanishes when the right member is in
$(\nabla_t M_A^\dual) \otimes ds$. We thus obtain a second pairing
\[
\frac{\{\bv \in M \otimes \calS: \nabla_t \bv \in M \otimes \calS^{t-} \otimes dx\}}
{M \otimes \calS^{t-}} \times ((R^1 f_* M_A^\dual) \otimes ds) \to K(-2)
\]
which is again nondegenerate on the left.
If we restrict further on the left to the classes of those $\bv \in
M \otimes \calS$ with
$\nabla_s \bv = 0$, then the pairing vanishes when the right
member is in $\calR^s_K \otimes_A \nabla_s (R^1 f_* M_A^\dual)$. 
We thus obtain a third pairing
\begin{equation} \label{eq:feqpair}
\frac{\{\bv \in M \otimes \calS: \nabla_s \bv = 0,
\nabla_t \bv \in M \otimes \calS^{t-} \otimes dx\}}
{\{\bv \in M \otimes \calS^{t-}: \nabla_s \bv = 0\}} 
\times H^1_{\loc,s=0}(R^1 f_* M_A^\dual) \to K(-2)
\end{equation}
which is again nondegenerate on the left. Combining this pairing
with the inclusion \eqref{eq:feqincl} yields the desired result.
\end{proof}

\begin{remark} \label{rem:degen}
In an earlier version of this paper, we attempted to construct
the embedding of Theorem~\ref{thm:degen} more directly, rather than
deduce it from Poincar\'e duality. This ran into trouble because the
convergence regions defining the rings $\calR^x_A$, for $A$ a
localization of $K \langle s \rangle^\dagger$, and
$\calR^s_{K \langle x \rangle^\dagger}$ do not share any common
territory. Thus we must avoid $R^1 f_! M_A$ and work with its dual instead,
which can be computed in the context of dagger algebras.
\end{remark}

\subsection{More degeneration in families}
\label{subsec:moredegen}

We continue to consider the situation of the previous section,
particularly in the case when the injection of Theorem~\ref{thm:degen}
is actually a bijection. We retain all notation from the previous section.
In addition, for $K'$ a finite extension of $K$, we write $M'
= M \otimes_K K'$, and for $\mu$ in the ring of integers $\gotho'$ of $K'$
fixed by a power of $\sigma_K$, we write $M_\mu = M'/(s-\mu)M'$.

\begin{lemma} \label{lem:factor}
Let $W$ be an $n \times n$ invertible matrix over $\calR^s_K$.
Then there exist $n \times n$ matrices $U$ and $V$ such that
$U$ is invertible over $\calR^{s,+}_K$, $V$ is invertible over
a localization of $K \langle s^{-1} \rangle^\dagger$
and $W = UV$.
\end{lemma}
\begin{proof}
Choose $r>0$ such that $w_r(W)$ is defined. (Note: in this argument, applying
$w_r$
to a matrix means taking its minimum over entries, rather than computing
like an operator norm.)
Choose a matrix $X$ over $K[s, s^{-1}]$ with nonzero determinant
such that $w_r(X - W^{-1}) > -w_r(W)$; then $w_r(WX- I) > 0$.
By \cite[Proposition~6.5]{bib:me7}, we can factor $WX$
as $YZ$, with $Y$ invertible over $\calR^{s,+}_K$ and $Z$ invertible over
$K \langle s^{-1} \rangle^\dagger$.
Since $\det(X) \in K[s,s^{-1}]$, $\det(X)$ is a unit in
some localization $A$ of $K \langle s^{-1} \rangle^\dagger$.
Thus we may put $U = Y$ and $V = ZX^{-1}$.
\end{proof}

\begin{lemma} \label{lem:excision}
Let $M$ be a $(\sigma, \nabla)$-module over 
$K \langle s,x \rangle^\dagger$ such that $H^0(M'_\mu) = 0$ for all
$K'$ and $\mu$.
Let $A \subseteq B$ be localizations of $K \langle s \rangle^\dagger$,
and suppose that $\bv \in M \otimes B \langle x \rangle^\dagger$ satisfies
$\nabla_t \bv \in M \otimes A \langle x \rangle^\dagger
\otimes dx$. Then $\bv \in M \otimes A \langle x \rangle^\dagger$.
\end{lemma}
\begin{proof}
Put $C = B \cap \calR^{s,+}_K$; it suffices to check that if $s$ is not 
invertible in $A$, then $\bv \in M \otimes C \langle x \rangle^\dagger$.
Namely, once this is done, we can repeat the argument after translating
(and enlarging $K$ as needed) to deduce that $\bv \in M \otimes A
\langle x \rangle^\dagger$.

By Proposition~\ref{prop:reldwork},
we can find a basis $\bv_1, \dots, \bv_n$ of $M \otimes 
\calR^{s,+}_{K \langle x \rangle^\dagger}$ such that $\nabla_s \bv_i = 0$
for $i=1, \dots, n$. From the integrability of $\nabla$, if we write
$\nabla_t \bv_j = \sum_i D_{ij} \bv_i \otimes dx$, then
$D_{ij} \in K \langle x \rangle^\dagger$ for all $i,j$.

Write $\bv = \sum_i a_i \bv_i$ with $a_i \in \calR^{s}_{K
\langle x \rangle^\dagger}$ (that is possible because
$\bv \in M \otimes B \langle x \rangle^\dagger$ and
$B \langle x \rangle^\dagger \subset \calR^s_{K \langle x \rangle^\dagger}$), 
write formally $a_i = \sum_l b_{il} s^l$
with each $b_{il} \in K \langle x \rangle^\dagger$,
and put $\bw_l = \sum_i b_{il} \bv_i$.
Then the series $\sum_l s^l \nabla_t \bw_l$ converges (in the
fringe topology of $M \otimes A \langle x \rangle^\dagger$) to $\nabla_t \bv$,
and the fact that $\nabla_t \bv \in M \otimes \calR^{s,+}_{K \langle x
\rangle^\dagger} \otimes dx$ implies that $\nabla_t \bw_l = 0$
for $l < 0$.

However, the $K \langle x \rangle^\dagger$-span of the $\bv_i$ is
a $(\sigma, \nabla)$-module isomorphic to $M/sM$, which by assumption
has no horizontal sections. Hence $\nabla_t \bw_l = 0$ if and only if
$\bw_l = 0$. We conclude that $\bw_l = 0$ for $l < 0$, and so
\[
\bv \in M \otimes (A \langle x \rangle^\dagger \cap
\calR^{s,+}_{K \langle x \rangle^\dagger}) = M \otimes
C \langle x \rangle^\dagger,
\]
which as noted above suffices to yield the desired result.
\end{proof}

\begin{prop} \label{prop:free}
Let $M$ be a free $(\sigma, \nabla)$-module over $K\langle s,x
\rangle^\dagger$, and let $f: K \langle s \rangle^\dagger \to
K \langle s,x \rangle^\dagger$ denote the canonical inclusion.
Suppose that for some nonnegative integer $m$,
\[
\dim_{K'} H^0(M_\mu)= \dim_{K'} H^0(M^\dual_\mu) = 0,
\qquad \dim_{K'} H^1(M_\mu)= \dim_{K'} H^1(M^\dual_\mu) = m
\]
for all $K'$ and $\mu$. Then $R^1 f_* M$, $R^1 f_* M^\dual$,
$R^1 f_! M$, $R^1 f_! M^\dual$
are free of rank $m$ over $K \langle s \rangle^\dagger$.
\end{prop}
\begin{proof}
It is enough to show that $R^1 f_* M$ and $R^1 f_! M$ are
locally free of rank $m$ over $K \langle s \rangle^\dagger$,
or likewise after replacing $K$ by a finite extension;
by translation, it suffices to check in a neighborhood of $s=0$.
Let $A$ be a localization of $K \langle s \rangle^\dagger$ such that
the conclusion of Theorem~\ref{thm:pushfwd} holds;
we may as well assume that $s$ is invertible in $A$, else we are already done.

We first treat $R^1 f_* M$ and $R^1 f_* M^\dual$. 
Under our hypothesis, the pairing
\eqref{eq:canfeq} must be perfect, as must be \eqref{eq:feqpair}.
In fact, the pairing
\[
[\cdot, \cdot]: \frac{\{\bv \in M \otimes \calS^{s+}: \nabla_s \bv = 0,
\nabla_t \bv \in M \otimes \calS^{s+,t-} \otimes dx\}}
{\{\bv \in M \otimes \calS^{s+,t-}: \nabla_s \bv = 0\}} 
\times H^1_{\loc,s=0}(R^1 f_* M_A^\dual) \to K(-2)
\]
is also perfect.

Choose a basis $\bv_1, \dots, \bv_m$ of $H^1_c(M/sM)$. Then the map
$h: (R^1 f_* M_A^\dual) \otimes_A \calR^s_K \to (\calR^s_K)^m$ defined by
\[
\bw \mapsto \left(\sum_{l \in \ZZ} s^l [g(\bv_i), s^{-l} \bw \otimes ds]\right)_{1 \leq i \leq m}
\]
is an isomorphism of $\calR^s_K$-modules. By
Lemma~\ref{lem:factor}, for some localization $B$ of $A$,
we can find elements $\bw_1, \dots, \bw_m$ of $R^1 f_* M_B^\dual$
such that the images $h(\bw_1), \dots, h(\bw_m)$ lie in
$(\calR^{s,+}_K)^m$ and generate $(\calR^{s,+}_K)^m$ 
 over $\calR^{s,+}_K$.
Choose $\bx_j \in M_B^\dual \otimes dx$ whose image in $R^1 f_* M_B^\dual$
is equal to $\bw_j$; then $[g(\bv_i), s^{-l} \bx_j \otimes ds] = 0$ 
for $l <0$, so each $\bx_j$
lies in
\[
M^\dual \otimes (B \langle x \rangle^\dagger \cap \calS^{s+}) \otimes dx
= M^\dual \otimes C \langle x \rangle^\dagger \otimes dx,
\]
where we write $C = B \cap \calR^{s,+}_K$; this is a localization of
$K \langle s \rangle^\dagger$ in which $s$ is not invertible.

Now given any $\bw \in R^1 f_* M^\dual_C$ which is the image of some
$\bx \in M_C^\dual \otimes dx$, we can uniquely write
$\bw = \sum_j b_j \bw_j$ with $b_j \in B$. On the other hand,
in the equality $h(\bw) = \sum_j b_j h(\bw_j)$, $h(\bw)$ belongs to
$(\calR^{s,+}_K)^m$, and the $h(\bw_j)$ generate $(\calR^{s,+}_K)^m$
freely over $\calR^{s,+}_K$. Therefore each $b_j$ in fact belongs to
$B \cap \calR^{s,+}_K = C$.

That is, given $\bx \in M^\dual_C \otimes dx$, $\bx - \sum_j b_j \bx_j$
is an element of $M^\dual_C \otimes dx$ which vanishes in $R^1 f_* M^\dual_B$.
By Lemma~\ref{lem:excision} (and using the hypothesis on the vanishing
of $H^0$), $\bx - \sum_j b_j \bx_j$ already vanishes
in $R^1 f_* M^\dual_C$. Hence $R^1 f_* M^\dual_C$ is freely generated by the
$\bw_j$. As noted above, this suffices to imply that $R^1 f_* M^\dual$ is free
of rank $m$; by the same argument, $R^1 f_* M$ is free of rank $m$.

We next consider $R^1 f_! M$. Choose a basis of $R^1 f_* M^\dual$, and let
$\bv_1, \dots, \bv_m$ be the dual basis of $R^1 f_! M_A$. Then under the
residue pairing
\[
\frac{M \otimes \calR^t_A}{M \otimes A \langle x \rangle^\dagger} 
\times (M^\dual \otimes dx)
\to A,
\]
each $\bv_i$ always pairs into $K \langle s \rangle^\dagger$. 
Hence $\bv_i$ is represented by an element of $M \otimes \calR^t_{K \langle
s \rangle^\dagger}$, and so belongs to $R^1 f_! M$. Similarly, given any
element of $R^1 f_! M$, we can write it uniquely as an $A$-linear
combination of the $\bv_i$ and then observe that the coefficients actually
lie in $K \langle s \rangle^\dagger$.
 Thus $R^1 f_! M$ (and likewise $R^1 f_! M^\dual$)
is also free of rank $m$, as desired.
\end{proof}

\section{A $p$-adic Fourier transform}
\label{sec:fourier}

In this chapter, we construct a $p$-adic version of the geometric Fourier
transform in $\ell$-adic cohomology, as introduced by Deligne and employed
by Laumon to give an alternate derivation of the Weil II theorem. Our point
of view will be from the theory of arithmetic $\calD^\dagger$-modules, but for 
the simple-minded manipulations we have in mind, we do not
need any of the rather substantial theory currently available about
such objects. In fact, our first act once we have defined our Fourier
transform will be to relate it to a geometric construction which
makes no reference to $\calD^\dagger$-modules. Later in the chapter, 
we give a version of the Grothendieck-Ogg-Shafarevich formula that
helps govern that geometric situation.

Nothing in this chapter is original, though we have worked things
out explicitly to illustrate the concrete nature of the construction.
The $p$-adic Fourier transform was introduced by
Mebkhout \cite{bib:meb1}, who also first proposed 
imitating Laumon's proof of Weil II in $p$-adic cohomology.
The Fourier transform was
generalized by Berthelot and studied more closely by Huyghe
\cite{bib:huyg1}. In particular, the coincidence with a 
geometric Fourier transform under suitable
conditions is due to Huyghe \cite{bib:huyg2}.
(Beware that our normalization is cosmetically different from Huyghe's:
we work with $\pi^{-n} \frac{d}{dx^n}$
in lieu of $\frac{1}{n!} \frac{d}{dx^n}$.)
The Grothendieck-Ogg-Shafarevich formula is cobbled together from
work of various authors; see Section~\ref{subsec:euler}.

Throughout this chapter, we assume that $K$ contains a primitive
$p$-th root of unity. It is equivalent to assume that $K$ contains
a $(p-1)$-st root of $-p$, which we will call $\pi$.

\subsection{The ring $\calD^\dagger$}
\label{subsec:thering}

We construct the ring $\calD^\dagger$ as a ring of ``overconvergent
differential operators'' on $K \langle x \rangle^\dagger$. We
then show that $(\sigma, \nabla)$-modules are in fact modules over
this ring. This will allow us to define a Fourier transform
by the usual trick of interchanging multiplication and differentiation.

\begin{lemma} \label{lem:factl}
For any positive integer $n$, $n/(p-1) \geq v_p(n!) \geq n/(p-1) - 
\lceil \log_p (n+1) \rceil$.
\end{lemma}
\begin{proof}
Use the formula
\[
v_p(n!) = \sum_{i=1}^{\infty} \lfloor np^{-i} \rfloor
\]
to write
\[
\frac{n}{p-1} - v_p(n!) = \sum_{i=1}^\infty  np^{-i} - \lfloor np^{-i} \rfloor.
\]
Every summand is nonnegative, yielding $n/(p-1) \geq v_p(n!)$. On the 
other hand, if $m = \lfloor \log_p n \rfloor$, then each summand 
with $i \leq m$
is bounded above by $1 - p^{-i}$, while each summand
 with $i > m$ is bounded above
by $np^{-i} \leq (p^{m+1} - 1)p^{-i}$. Since
\[
\sum_{i=1}^m (1-p^{-i}) + \sum_{i=m+1}^\infty (p^{m+1}-1)p^{-i} = m+1,
\]
we have $n/(p-1) - v_p(n!) \leq m+1 = \lceil \log_p(n+1) \rceil$ as desired.
\end{proof}

Let $R$ be the noncommutative polynomial ring in the variables
$x$ and $\del$ over $K$, modulo the two-sided ideal generated by
$\del x - x \del - \pi^{-1}$.
\begin{lemma} \label{lem:noncom}
In $R$, one has
\[
\del^n x^m = \sum_i \frac{n!m!}{\pi^i i!(n-i)!(m-i)!} x^{m-i} \del^{n-i}.
\]
\end{lemma}
\begin{proof}
The equality is evident if $m=0$ or $n=0$, and the case where $m=1$ or
$n=1$ is also easily checked by induction on the other variable.
The general case follows
by induction on $m+n$:
\begin{align*}
  \del^n x^m &= (x\del^n + n \pi^{-1} \del^{n-1}) x^{m-1} \\
&= x \sum_i \frac{n!(m-1)!}{\pi^i i!(n-i)!(m-1-i)!} x^{m-1-i} \del^{n-i}\\
&\mbox{} \qquad +  n \pi^{-1} \sum_i 
\frac{(n-1)!(m-1)!}{\pi^i i!(n-1-i)!(m-1-i)!} x^{m-1-i} \del^{n-1-i} \\
&= \sum_i (n(m-i) + ni)
\frac{(n-1)!(m-1)!}{\pi^i i!(n-i)!(m-i)!} x^{m-i} \del^{n-i} \\
&= \sum_i \frac{n!m!}{\pi^i i!(n-i)!(m-i)!} x^{m-i} \del^{n-i}.
\end{align*}
\end{proof}

Let $\calD^{x,\dagger}$ 
be the set of formal power series $\sum_{i,j=0}^\infty a_{ij}
x^i \del^j$ such that $\liminf_{i,j} \{v_p(a_{ij})/(i+j)\} > 0$;
we drop the series parameter $x$ when it is understood.
(We will not define $\calD$ unadorned, but the dagger will help remind
us of the ``overconvergent'' nature of the definition.)

We use Lemma~\ref{lem:noncom} to show that multiplication of
power series in $\calD^\dagger$ is well-defined.
\begin{prop} \label{prop:mult}
Let $a = \sum_{i,j=0}^\infty a_{ij} x^i \del^j$ and
$b = \sum_{k,l=0}^\infty b_{kl} x^k \del^l$ be elements of $\calD^\dagger$.
Then the series
\[
c_{mn} = \sum_{i,j,s} a_{ij} b_{(m+s-i)(n+s-j)}
\frac{j!(m+s-i)!}{\pi^s s!(j-s)!(m-i)!}
\]
converges in $K$ for each $m,n$,
and the series $c = \sum_{m,n=0}^\infty c_{mn} x^m \del^n$ belongs to 
$\calD^\dagger$.
\end{prop}
\begin{proof}
By the definition of $\calD^\dagger$, there exist constants $e,f$ such that
$v_p(a_{ij}) \geq e(i+j)-f$ and $v_p(b_{kl}) \geq e(k+l) - f$  for all
$i,j,k,l$.
By Lemma~\ref{lem:factl}, the valuation of the summand for any
given $i,j,s$ is at least
\begin{align*}
& v_p(a_{ij}) + v_p(b_{(m+s-i)(n+s-j)}) + 
\frac{j+(m+s-i) - 2s-(j-s)-(m-i)}{p-1} \\
& \qquad - \lceil \log_p (j+1) \rceil - \lceil \log_p (m+s-i+1)\rceil \\
&\geq e(i+j+(m+s-i)+(n+s-j)) - 2f-2 - \log_p (j+1) - \log_p (m+s-i+1) \\
&= e(m+s) - \log_p (m+s-i+1) + e(n+s) - \log_p (j+1) - 2f-2 \\
&\geq e(m+s) - \log_p (m+s+1) + e(n+s) - \log_p (n+s+1) - 2f-2.
\end{align*}
This expression tends to infinity as $s \to \infty$, so the
sum converges. More precisely,
for any $g$ with $0<g<e$, there exists $h$ such that
$ex - \log_p (x+1) \geq gx - h$ for all $x \geq 1$.
Then
\begin{align*}
v_p(c_{mn}) &\geq \min_{s \geq 0} \{ g(m+s) - h + g(n + s) - h - 2f-2\} \\
&= g(m+n) - 2h - 2f-2.
\end{align*}
Thus $\sum_{m,n=0}^\infty c_{mn} x^m \del^n \in \calD^\dagger$, as desired.
\end{proof}
\begin{remark} \label{rem:switch}
By the same argument, each element of $\calD^{\dagger}$
can be written as $\sum a_{ij} \delta^j x^i$, again with
$\liminf v_p(a_{ij})/(i+j) > 0$. This will come up again when we define
the Fourier transform in the next section.
\end{remark}

In the notation of Proposition~\ref{prop:mult},
we define a multiplication operation on $\calD^\dagger$ by setting
$ab = c$. On the subset $R$ of $\calD^\dagger$, this operation coincides with
the multiplication in $R$ by Lemma~\ref{lem:noncom}. 
The usual ring axioms can thus be verified by approximating elements
of $\calD^\dagger$ with elements of $R$. (More precisely, one can give
$\calD^\dagger$ a ``fringe topology'' like that of $W_n$, under
which $R$ is visibly dense, and the proof of Proposition~\ref{prop:mult}
shows that multiplication is continuous.)
The upshot is that $\calD^\dagger$ forms a (noncommutative)
ring under series multiplication.

We note in passing that the ring $\calD^\dagger$ coincides with the
``overconvergent Weyl algebra''
$A_1(K)^\dagger$ considered in \cite{bib:huyg3}. There it is shown
that the category of coherent left $A_1(K)^\dagger$-modules coincides
with the category of coherent left modules for the sheaf
$\calD^{\dagger}_{\AAA^1, \QQ}(\infty)$ of noncommutative rings,
constructed by Berthelot using divided power envelopes.

\subsection{$\calD^{\dagger}$-modules and the Fourier transform}

To construct the Fourier transform of a $(\sigma, \nabla)$-module
over $K \langle x \rangle^\dagger$, we
must first convert it into a (left) $\calD^\dagger$-module, in which $\del$
acts so that $\del \bv \otimes dx = \pi^{-1} \nabla \bv $. 
However, it is not immediately
obvious that this action makes sense for \emph{power series} in $\del$,
so we must verify this first. 

The conclusion of the following lemma is essentially part of the definition of
an isocrystal without Frobenius structure; the lemma says that this
part of the definition is superfluous in the presence of Frobenius.
\begin{lemma} \label{lem:converge}
Let $M$ be a $(\sigma, \nabla)$-module over $K \langle x \rangle^\dagger$,
and define $D: M \to M$ by $\nabla \bv = D\bv\otimes dx$.
Then for any sequence $\{\bv_i\}_{i=0}^\infty$ of elements of $M$
convergent to zero under the fringe topology, and any $\epsilon \in
K \langle x \rangle^\dagger$
with $|\epsilon| < 1$, the double sequence
$\{\epsilon^j \pi^{-j} D^j \bv_i\}_{i,j=0}^\infty$ converges to zero
under the fringe topology of $M$.
\end{lemma}
\begin{proof}
Choose an integer $e$ large enough so that
\[
1 > |\epsilon/\pi|^{p^e} |(p^e)!| =
|\epsilon|^{p^e} p^{1/(p-1)}.
\]
Choose a basis $\be_1, \dots, \be_n$ of $M$, and define the matrix $N$ by
\[
D\be_j = \sum_i N_{ij} \be_i.
\]
Put $u = \frac{dx^\sigma}{dx}$ and define the matrix $N^{(m)}$ by
\[
N^{(m)} = u u^\sigma \cdots u^{\sigma^{m-1}} N^{\sigma^m};
\]
then we have $DF^m \be_j = \sum_i N^{(m)}_{ij} F^m \be_i$.
Since $|u| < 1$, we may choose $m$ large enough so that
$|N^{(m)}| < |(p^e)!|$.

Write $\bv_i = \sum_l c_{il} F^m \be_l$; then as $i \to \infty$,
the $c_{il}$ converge to
zero within $T_{1,\rho}$ for some $\rho > 1$. 
Let $|\cdot|_\rho$ denote the spectral norm
on $T_{1,\rho}$ (i.e., the Gauss norm). 
By choosing $\rho$
sufficiently close to 1, we can ensure that 
$|\epsilon/\pi|_\rho^{p^e} |(p^e)!| < 1$ and that
$N^{(m)}/(p^e)!$ has entries
in $T^{\inte}_{1,\rho}$.
Writing $D^{p^e} \bv_i = \sum_l 
d_{il} F^m \be_l$, we obtain
\begin{align*}
|d_{il}|_{\rho} &\leq \max\{|\frac{d^{p^e}}{dt^{p^e}} 
c_{il}|_{\rho}, 
|N^{(m)}|_\rho \max_l \{|c_{il}|_\rho\}\} \\
&\leq |(p^e)!| \max_l \{|c_{il}|_\rho\}.
\end{align*}
If we write $D^j \bv_i = \sum_l f_{ijl} F^m \be_l$, we may deduce that
\begin{align*}
|\epsilon^j \pi^{-j} d_{ijl}|_\rho &\leq
|\epsilon^j \pi^{-j} (p^e)!^{\lfloor j/p^e \rfloor}|_\rho
\max_l \{|c_{il}|_\rho\} \\
&\leq |\epsilon/\pi|_\rho^{j - p^e \lfloor j/p^e \rfloor}
|(\epsilon/\pi)^{p^e} (p^e)!|_\rho^{\lfloor j/p^e \rfloor}
\max_l \{|c_{il}|_\rho\} \\
&\leq \max\{1, |\epsilon/\pi|_\rho^{p^e-1}\}
|(\epsilon/\pi)^{p^e} (p^e)!|_\rho^{\lfloor j/p^e \rfloor}
\max_l \{|c_{il}|_\rho\}
\end{align*}
which converges to zero as $i+j \to \infty$. This yields the desired
convergence.
\end{proof}

\begin{cor} \label{cor:double series}
Let $M$ be a $(\sigma, \nabla)$-module over $K \langle x \rangle^\dagger$,
and define $D: M \to M$ by $\nabla\bv = D\bv\otimes dx$.
For any $\sum_{i,j} a_{ij} \del^j x^i \in \calD^\dagger$ and
any $\bv \in M$, the double series $\sum_{i,j} a_{ij}
(\pi^{-1} D)^j x^i \bv$
converges in $M$.
\end{cor}
\begin{proof}
We can find $\delta$ in a finite extension of $K$ with $\delta < 1$
such that $|a_{ij}|<|\delta|^{i+j}$ for all but finitely many pairs $i,j$.
Now apply Lemma~\ref{lem:converge} with $\epsilon = \delta$
and $\bv_i = (\delta x)^i \bv$.
\end{proof}
By Corollary~\ref{cor:double series}, any $(\sigma, \nabla)$-module
over $K \langle x \rangle^\dagger$ can be given the structure
of a left $\calD^\dagger$-module in which
\[
\del\bv\otimes dx = \pi^{-1} \nabla\bv.
\] 

The ring $\calD^\dagger$ admits an automorphism $\rho$ sending $x$ to
$\del$ and $\del$ to $-x$. More explicitly, the formula for
$\rho$ can be read off from Lemma~\ref{lem:noncom}:
\[
\rho\left( \sum_{i,j} c_{ij} x^i \del^j \right)
 =
\sum_{i,j} \left( \sum_k (-1)^{j+k}
\frac{(i+k)!(j+k)!}{\pi^k k! i!j!} c_{(i+k)(j+k)} \right)  x^i \del^{j}.
\]
Thus given a left $\calD^\dagger$-module $M$, we get
a new left 
$\calD^\dagger$-module $\calD^\dagger \otimes_{\rho} M$; we call this the
\emph{Fourier transform} of $M$ and denote it by $\widehat{M}$.

\begin{remark} \label{rem:buildinF}
The Frobenius structure on the Fourier transform of a $(\sigma, \nabla)$-module
will be obtained using the geometric Fourier transform in the next
section. However, it is worth sketching an alternate approach:
 one can enlarge $\calD^\dagger$ to include an element $F$ satisfying
the relations
\[
F x = x^\sigma F \qquad \mbox{and} \qquad
\del F = \frac{dx^\sigma}{dx} F \del.
\]
The automorphism $\rho$ can be extended to the larger ring as follows.
Define $c_i \in K$ by the formal identity
\[
\sum_{i=0}^\infty c_i x^i = \exp(- \pi x + \pi x^\sigma);
\]
then in fact $\sum c_i x^i \in K \langle x \rangle^\dagger$
(an observation of Dwork),
and the extension of $\rho$ satisfies
\[
\rho(F) = \left( \sum_{i=0}^\infty c_i \del^i x^i \right) 
\frac{dx^\sigma}{dx} F
\]
and $\rho(\rho(F)) = qF$. 
\end{remark}

\begin{remark}
In fact, the ring constructed in Remark~\ref{rem:buildinF} has
an $F$ corresponding to each Frobenius lift $\sigma$, giving a nice
interpretation of the fact that the category of $(\sigma, \nabla)$-modules
on $K \langle x \rangle^\dagger$ does not depend on the choice of 
the Frobenius lift $\sigma$.
\end{remark}

\subsection{The (na\"\i ve) geometric Fourier transform}

The description of the Fourier transform given above is concise
and elegant, but not particularly amenable to analysis of the sort
we wish to carry out. For this, we need a more explicit description;
we get this description at the expense of restricting $M$.

The Dwork isocrystal (or Artin-Schreier isocrystal) on the $x$-line
$\AAA^1$ is an overconvergent $F$-isocrystal of rank one,
defined as follows. Associating to $\AAA^1$ the dagger
algebra $K \langle x \rangle^\dagger$ with its standard Frobenius,
we define a $(\sigma,
\nabla)$-module $\calL$ of rank one over $K \langle x \rangle^\dagger$ by
giving a single generator $\be$ and the Frobenius and connection
actions
\[
F\be = \exp(\pi x - \pi x^q)\be, \qquad
\nabla \be = \pi \be \otimes dx.
\]
This isocrystal becomes trivial after adjoining $u$ such that $u^p - u = x$.
(This implies that its $p$-th tensor power is already trivial on the
$x$-line.)

For any dagger algebra $A$ and any $f \in A^{\inte}$, we identify
$f$ with the map $K \langle x \rangle^\dagger \to A$ 
mapping $x$ to $f$, and write
$\calL_f$ for $f^* \calL$. Note that $\calL_{f+g} = \calL_f \otimes \calL_g$
and that the isomorphism class of $\calL_f$ depends only on $f$ modulo $\pi$.

Let $M$ be a $(\sigma, \nabla)$-module over $K \langle x \rangle^\dagger$,
and let $f: K \langle s \rangle^\dagger
\to K \langle s,x \rangle^\dagger$ and
$g: K \langle x \rangle^\dagger \to K \langle s,x \rangle^\dagger$
be the canonical embeddings.
Then $g^* M$ and $\calL_{sx}$ are $(\sigma, \nabla)$-modules over
$K \langle s,x \rangle^\dagger$, as is
\[
N = g^* M \otimes_{K \langle s,x \rangle^\dagger} \calL_{sx}.
\]
We can decompose $\Omega^1_{K \langle s,x\rangle^\dagger}$ into
two rank one submodules, generated by $ds$ and $dx$. Let
$\nabla_s$ and $\nabla_x$ be the components of the connection
on $N$ mapping to these two submodules.

We define the \emph{(na\"\i ve) geometric Fourier transform}
of $M$ as $\widehat{M}_{\geom} = \coker \nabla_x$; this is a
$\calD^{s,\dagger}$-module equipped with a $\sigma$-linear Frobenius map,
but is not necessarily locally free.
The nomenclature
is justified by the following fact.

\begin{prop} \label{prop:notnaive}
Let $M$ be a $(\sigma, \nabla)$-module over $K \langle x \rangle^\dagger$.
Then there is a canonical isomorphism
$\widehat{M} \to \widehat{M}_{\geom}$ of $\calD^\dagger$-modules.
\end{prop}
\begin{proof}
The map in question is defined as follows.
We identify $\widehat{M}$ with $M$ as sets (or even as $K$-vector spaces)
by identifying $\bv \in M$
with $1 \otimes \bv$. We then identify $M$ with a subset of $g^* M$ via $g$,
and in turn identify $g^* M$ with $N$ by identifying $\bw \in g^* M$
with $\bw \otimes \be$ (where $\be$ is the distinguished generator
used in the definition of $\calL$, or more precisely, its image in
$\calL_{sx}$). The desired map is now constructed by tracing through
these identifications, then composing with the map $N \to \coker \nabla_x$
induced by $\bv \mapsto \bv \otimes dx$.

We check that this map is surjective.
Let $D: M \to M$ be the map with $\nabla\bv = D\bv \otimes dx$.
Given $\bv \in N$, we may write
$\bv = \sum_i \bv_i s^i$ for some
$\bv_i \in M$. 
In this representation, the sequence $\eta^i \bv_i$ must converge to zero 
(under the fringe topology of $M$) for
some $\eta$ in a finite extension of $K$ with $|\eta| > 1$.
Apply Lemma~\ref{lem:converge} to the sequence $\{\eta^i 
\bv_i\}$ with $\epsilon = \eta^{-1}$ to deduce that
$\eta^{i-j} \pi^{-j} D^j  \bv_{i}$ converges to zero.
The same holds
if we replace $i$ by $i+j+1$ (since this restricts the double sequence), 
that is, $\eta^{i+1} \pi^{-j} D^j  \bv_{i+j+1}$ converges to zero.
Thus the series
\[
\sum_{i=0}^\infty s^i
\sum_{j=0}^\infty (-1)^j \pi^{-(j+1)} D^j
\bv_{i+j+1}
\]
converges in $M$ to a limit $\bw$. 
Because of the way we identified $M$ within $N$,
we have
\[
\nabla_x D^j \bv_{i+j+1} = (D^{j+1} \bv_{i+j+1} + \pi s D^j \bv_{i+j+1}) 
\otimes dx,
\]
from which it follow that 
\[
\nabla_x \bw = (\bv - \bv_0)\otimes dx
+ \sum_{j=0}^\infty (-1)^j \pi^{-j-1} D^{j+1} \bv_{j+1} \otimes dx.
\]
We conclude that
every element of $\coker \nabla_x$
is represented by an element of the form $\bv_0\otimes dx$ with $\bv_0 \in M$,
so the map $M \to \coker \nabla_x$ is surjective.

We next check injectivity. Suppose $\bv \in M$ becomes
zero in $\coker \nabla_x$; that means there exists
$\bw = \sum_i \bw_i s^i \in N$
such that $\nabla_x \bw = \bv\otimes dx$.
Comparing coefficients of powers of $s$, we have
$D \bw_0 = \bv$ and
$\pi \bw_i + D \bw_{i+1} = 0$ for $i \geq 0$,
implying $(-\pi)^{-i} D^{i+1} \bw_i = \bv$.
On the other hand, we can choose $\eta$ in a finite extension of $K$
with $|\eta| > 1$ such that the sequence
$\{\eta^i \bw_i\}$ converges to zero in $M$.
Then by Lemma~\ref{lem:converge} applied with $\epsilon = \eta^{-1}$
and $\bv_i = \eta^i \bw_i$, $(-\pi)^{-i} D^i \bw_i$ converges
to zero in $M$, a contradiction unless $\bv=0$. Thus the map
is injective.
\end{proof}

For our purposes, the main significance of this result is
the following, which we state in the notation of
Section~\ref{subsec:moredegen}.
\begin{prop} \label{prop:irred}
Let $M$ be a $(\sigma, \nabla)$-module over $K \langle x \rangle^\dagger$.
Suppose there exists a nonnegative integer $d$ such that
that for each $K'$ and $\mu$,
\begin{gather*}
\dim_{K'} H^0(M \otimes \calL_{\mu x}) = 
\dim_{K'} H^0(M^\dual \otimes \calL_{\mu x}) = 0, \\
\dim_{K'} H^1(M \otimes \calL_{\mu x}) =
\dim_{K'} H^1(M^\dual \otimes \calL_{\mu x}) = d.
\end{gather*}
Then $\widehat{M}_{\geom}$ is
a $(\sigma, \nabla)$-module of rank $d$
over $K \langle s \rangle^\dagger$.
If in addition $M$ is (absolutely) irreducible
as a $(\sigma, \nabla)$-module,
then $\widehat{M} \cong \widehat{M}_{\geom}$ is (absolutely) irreducible
as a $(\sigma, \nabla)$-module.
\end{prop}
\begin{proof}
The first assertion follows immediately from Proposition~\ref{prop:free},
so we focus on the second.
If $\widehat{M}$ is reducible as a $(\sigma, \nabla)$-module,
it has a Frobenius-stable
$\calD^\dagger$-submodule $\widehat{N}$ such that $\widehat{N}$
and $\widehat{M}/\widehat{N}$
are infinite dimensional $K$-vector spaces. Undoing the 
Fourier transform gives a Frobenius-stable
$\calD^\dagger$-submodule $N$ of $M$ such that
$N$ and $M/N$ are infinite dimensional $K$-vector spaces. 
But then $N$ is a $(\sigma, \nabla)$-submodule of $M$ such that
$N$ and $M/N$ are nontrivial, so $M$ is reducible.
Hence if $M$ is irreducible, then so is $\widehat{M}$; the same is
true with ``irreducible'' replaced by ``absolutely irreducible'' because the
construction of the $\calD^{\dagger}$-module
Fourier transform clearly commutes with extension of the base field.
\end{proof}

\subsection{An Euler characteristic formula}
\label{subsec:euler}

In order to apply the results of the previous section to a
$(\sigma, \nabla)$-module $M$ over $K \langle x \rangle^\dagger$, we need
to establish conditions under which the dimension of $H^1(M \otimes
\calL_{\mu x})$ does not depend on $\mu$. This requires
a formula for this dimension; in this section,
we establish such a formula
using some recent results in $p$-adic cohomology.

The $\ell$-adic analogue of the formula we seek
is the Grothendieck-Ogg-Shafarevich formula
\cite{bib:sga5} (see also \cite{bib:raynaud}), which relates the 
Euler-Poincar\'e
characteristic of a lisse sheaf on a curve to the local monodromy at
the missing points. Naturally, its $p$-adic analogue will also be given
in terms of local monodromy.

Let $C$
 be a smooth irreducible affine curve over $k$,
let $\overline{C}$ be the smooth
compactification of $k$, and let $\calE$ be an overconvergent $F$-isocrystal
on $C$. 
Let $x$ be a closed point of $\overline{C} \setminus C$,
let $E$ be
the fraction field of the completed local ring of $C$ at $x$, and let $F$
be a Galois extension of $E$ over which the local monodromy of $\calE$
at $x$ becomes unipotent (or more precisely, over which 
the module obtained from a
$(\sigma, \nabla)$-module corresponding to $\calE$ by tensoring up
to a Robba ring $\calR_x$ corresponding to $x$ becomes unipotent).
Let $f$ be the residue field degree of $F/E$,
and let $G$ be the Galois group $\Gal(F/E)$.
Define the Swan function on $G$
by the formula
\[
\Swan_{F/E}(g) = \begin{cases}
- f \inf_{x \in \gotho_F \setminus \{0\}} \{v_F(x^g/x - 1)\} & 
g \neq e \\
- \sum_{h \in G \setminus \{e\}} \Swan_{F/E}(h) & g = e.
\end{cases}
\]
Let $\chi: G \to K$ be the character of the representation of $G$
on the local monodromy of $\calE$ at $x$, i.e., on the horizontal sections
of $\calE$ over the extension of $\calR_x$ corresponding to $F$
(and the horizontal sections of the quotient of $\calE$
by the span of all horizontal sections, and so on).
Define the \emph{Swan conductor} of $\calE$ at $x$ as
\[
\Swan_x(\calE) = \frac{1}{|G|} \sum_{g \in G} \Swan_{F/E}(g) \chi(g);
\]
it turns out that this quantity is an integer (by a theorem of Artin),
and does not depend on the choice of $F$.

A more convenient recipe for computing the Swan conductor is the
following. (See \cite{bib:serre}
and/or \cite[Chapter~1]{bib:katzgauss} for more details.) 
For $i \geq 0$, let $G^i$ be the $i$-th ramification subgroup
in the upper numbering. Given an irreducible 
representation $\rho: G^0 \to \GL(V)$ with open kernel on a finite
dimensional $K$-vector space, the smallest number
$i$ (necessarily rational)
such that $G^i \subseteq \ker(\rho)$ is called the \emph{break} of $\rho$.
For $\rho$ the local monodromy representation of $\calE$ at $x$,
we can decompose $V = \oplus_{i \geq 0} V(i)$, where
$V$ is the direct sum of the irreducible subrepresentations of break
$i$; we then have the formula
\[
\Swan(\rho) = \sum_{i \geq 0} i \dim V(i).
\]

In terms of the Swan conductor, the desired analogue of the
Grothendieck-Ogg-Shafarevich formula is as follows.
\begin{theorem} \label{thm:gos}
Let $\calE$ be an overconvergent $F$-isocrystal on a smooth irreducible
affine curve $C$ over $k$, and write
\[
\chi(C/K,\calE) = \dim_K H^0_{\rig}(C/K,\calE) - \dim_K H^1_{\rig}(C/K,\calE).
\]
Then
\begin{equation} \label{eq:rigep}
\chi(C/K,\calE) = \chi(C/K, \calO_C) \rank(\calE)
- \sum_{x \in \overline{C} \setminus C} [\kappa(x):k] \Swan_x(\calE).
\end{equation}
\end{theorem}
\begin{proof}
The theorem  may be obtained at once by combining
the following two results.
\begin{enumerate}
\item[(a)] A theorem of Christol and Mebkhout 
\cite[Corollaire~5.0--12]{bib:cm} 
states that
\[
\chi(C/K,\calE) = 
\chi(C/K, \calO_C) \rank(\calE) 
- \sum_{x \in
\overline{C} \setminus C} [\kappa(x):k] \Irr_x(\calE),
\]
where $\Irr$ is the ``irregularity'' of $\calE$ at $x$
(a generalized form of a definition of Robba).

\item[(b)] 
A theorem of Crew \cite[Theorem~5.4]{bib:crew3},
Matsuda \cite[Theorem~8.6]{bib:matsuda},
and Tsuzuki \cite[Theorem~7.2.2]{bib:tsuswan} states that
\[
\Irr_x(\calE) = 
\Swan_x(\calE).
\]
\end{enumerate}
However, in keeping with the semi-expository style of this paper,
we sketch a proof obtained by combining the proofs of (a) and (b)
and then streamlining.

We start with some reductions.
First, note that there is no harm at any point in replacing $K$ and $k$
by finite
extensions. Second, note that it
suffices to prove the claim after shrinking $C$: it suffices to observe that
removing a single rational point increases both sides of
\eqref{eq:rigep} by $\rank(\calE)$ (since the irregularity is zero at the
removed point).
Third, note that it suffices to prove the claim
after pushing forward along a finite \'etale map $C \to D$. In particular,
by a technique of Abhyankar (included in Proposition~\ref{prop:belyi}),
we may reduce to the case $C = \AAA^1$, in which case \eqref{eq:rigep}
reduces to
\[
\chi(\AAA^1/K,\calE) = \rank(\calE)
- \Swan_\infty(\calE)
\]
or equivalently
\[
\chi(\GG_m/K,\calE) = - \Swan_\infty(\calE),
\]
where $\GG_m = \AAA^1 \setminus \{0\}$.

By \cite[Th\'eor\`eme~5.0--10]{bib:cm}, one can extend $\calE$ to
a rigid analytic vector bundle $\calF$ on $\PP^1_K$ 
equipped with a (possibly irregular) meromorphic connection
on $\AAA^1_K = \PP^1_K \setminus \{\infty\}$ which is holomorphic
away from $\infty$. (This is a local claim in the residue disc at
$\infty$, so one can prove it by reducing to the case of an irreducible
local monodromy representation; something like this indeed occurs in
\cite[Th\'eor\`eme~5.0--10]{bib:cm}, giving a reduction to
\cite[Th\'eor\`eme~4.2--7]{bib:cm}.)
By rigid-analytic GAGA, this bundle and connection are both algebraic,
and its algebraic and rigid-analytic Euler characteristics coincide.

Let $N$ be the module over $K[x,x^{-1}]$ 
corresponding to $\calF$ over $\GG_{m,K}$,
so that
$N \otimes K \langle x,x^{-1} \rangle^\dagger$ is the $(\sigma, \nabla)$-module
corresponding to $\calE$ over $k[x,x^{-1}]$. 
Write 
\begin{align*}
R_0 &= K \langle x^{-1} \rangle^\dagger/K[x^{-1}] =
\calR^{x}/x^{-1} \calR^{x,+} \\
R_\infty &= K \langle x \rangle^\dagger/K[x] =
\calR^{x^{-1}}/x \calR^{x^{-1},+};
\end{align*}
then there is a natural isomorphism
\[
K \langle x,x^{-1} \rangle^\dagger / K[x,x^{-1}] \cong R_0 \oplus R_\infty.
\]
Thus we have a commuting diagram with exact rows
\[
\xymatrix{
0 \ar[r] & N \ar^\nabla[d] \ar[r] & N \otimes K \langle x,x^{-1} \rangle^\dagger
\ar^\nabla[d] \ar[r] & N \otimes (R_0 \oplus R_\infty) 
\ar^\nabla[d] \ar[r] & 0 \\
0 \ar[r] & N \otimes \Omega^1 
 \ar[r] & N \otimes K \langle x,x^{-1} \rangle^\dagger \otimes 
\Omega^1
 \ar[r] & N \otimes (R_0 \oplus R_\infty) \otimes \Omega^1 
\ar[r] & 0,
}
\]
from which the snake lemma gives
\begin{equation} \label{eq:snakelemma}
\chi(\GG_{m}/K,\calE) - \chi(\GG_{m,K},\calF) =
\chi(N \otimes R_0) + \chi(N \otimes R_\infty).
\end{equation}

Following Crew, Matsuda, and Tsuzuki,
one now constructs
a local-to-global correspondence in the vein of
Katz's \cite[Main~Theorem~1.4.1]{bib:katzlocal}. 
This produces a unit-root overconvergent $F$-isocrystal $\calE'$ on
$\GG_m$ having the same local monodromy at $\infty$ as does $\calE$,
and having tame monodromy at 0. We can extend $\calE'$ to
a rigid analytic vector bundle $\calF'$ on $\PP^1_K$ 
with irregular connection on $\GG_{m,K}$ by using the \emph{same}
extension to the residue disc at $\infty$ as was taken for $\calF$,
and using an extension at 0 which makes the connection regular singular
(i.e., at worst simple poles).
We then have an analogue of \eqref{eq:snakelemma}:
\[
\chi(\GG_{m}/K,\calE') - \chi(\GG_{m,K},\calF') =
\chi(N' \otimes R_0) + \chi(N' \otimes R_\infty).
\]
However, $N' \otimes R_\infty \cong N \otimes R_\infty$ by construction,
and $\chi(N' \otimes R_0) = 0$ because the connection is regular singular
at 0.
Hence
\[
\chi(\GG_{m}/K,\calE') - \chi(\GG_{m,K},\calF') =
\chi(\GG_{m}/K,\calE) - \chi(\GG_{m,K},\calF).
\]
Since $\calE$ and $\calE'$ both have the same Swan conductor at $\infty$
and at $0$ (the ones at $0$ both being zero), 
it now suffices to make the following two observations.
\begin{enumerate}
\item[(i)] The formula \eqref{eq:rigep} holds with $\calE$ replaced by
$\calE'$. This can be checked in two ways paralleling two arguments
in \'etale cohomology: one may pass to crystalline cohomology and
use intersection theory, for which see \cite[\S 4]{bib:crew3};
or one may perform a direct calculation using Brauer induction,
for which see \cite[\S 8]{bib:tsuswan}. (In particular, the second
approach maintains our informal self-prohibition against using
crystalline methods.)
\item[(ii)] The equality $\chi(\GG_{m,K},\calF') = \chi(\GG_{m,K},\calF)$
holds. By complex-analytic GAGA (after choosing an embedding
$K \hookrightarrow \CC$), this may be deduced
from Mebkhout's Euler-Poincar\'e 
formula for a compact Riemann surface equipped with meromorphic connection
\cite[\SS 2.3]{bib:meb0}, as follows.
The formula asserts that
\[
\chi(\GG_{m,K}, \calF) = -\Irr_0(\calF) - \Irr_\infty(\calF)
\]
and similarly for $\calF'$, where $\Irr_0$ and $\Irr_\infty$ denote
the irregularity in the sense of Malgrange \cite{bib:malgrange}.
On one hand, $\Irr_0(\calF) = \Irr_0(\calF') = 0$ because
$\calF$ and $\calF'$ are both regular at $0$; on the other hand,
$\Irr_\infty(\calF) = \Irr_\infty(\calF')$ because at $\infty$,
$\calF$ and $\calF'$ are isomorphic rigid-analytically in a neighborhood, 
hence formally, hence complex-analytically in a neighborhood.
(Mebkhout attributes this Euler-Poincar\'e formula
to Deligne, but the citation of 
\cite{bib:del0} in \cite{bib:meb0} is garbled.)
\end{enumerate}
This yields the desired result.
\end{proof}

\begin{remark} \label{rem:gos}
It would be possible to give a more direct
proof of Theorem~\ref{thm:gos} (one not relying on complex analysis)
if one could establish Oort's conjecture
\cite[1.6]{bib:oort1} (see also \cite{bib:gm})
that for $k$ algebraically closed,
every finite cyclic cover between smooth proper curves over $k$
lifts to such a cover over $\gotho$. In that case, one can
establish an equivariant form of
Grothendieck-Ogg-Shafarevich (as in \cite{bib:raynaud}) by invoking the
$p$-adic monodromy theorem to reduce to
the case where $\calE$ has everywhere unipotent monodromy, then
computing explicitly as in \cite[Chapter~6]{bib:me8} to establish
the formula for one group element at a time. Note that in this situation, the 
Euler-Poincar\'e characteristics coincide in the rigid cohomological
and algebraic/analytic settings; it is conceivable that
inability to achieve this coincidence
could constitute an obstruction to the truth of Oort's conjecture,
or even to its weak form in which $\gotho$ may be chosen depending on
the cover.
\end{remark}

Using Theorem~\ref{thm:gos}, we obtain the following result.
\begin{prop} \label{prop:swan}
Let $M$ be a $(\sigma, \nabla)$-module over $K\langle x\rangle^\dagger$.
Then there exists an integer $N$ such that
for $P \in \gotho[x]$ a monic polynomial of degree $d$, with $d
> N$ and $d$ not divisible by $p$, 
we have
\begin{gather*}
\dim_K H^0_{\loc}(M \otimes \calL_P) = \dim_K H^1_{\loc}(M \otimes \calL_P)
= 0, \\
\dim_K H^0(M \otimes \calL_P) = 0, \qquad
\dim_K H^1(M \otimes \calL_{P}) = (d - 1) \rank(M).
\end{gather*}
\end{prop}
\begin{proof}
Let $\rho: G \to \GL(V)$ 
be the local monodromy representation of $M$ at infinity;
then the local monodromy representation of $M \otimes \calL_P$ at infinity
is equal to
$\rho \otimes \psi_P$, for $\psi_P$ a suitable nontrivial
character of the Galois
group $\Gal(L/k((t)))$, with $L = k((t))[u]/(u^p - u - P(t))$. 
Let $N$ be the largest break of $\rho$; this will turn out to be a good
choice. 

Suppose $d > N$.
Then $\rho \otimes \psi_P$ has all breaks equal to $d$,
because the subgroup of $G$ fixing $L$ is not contained in $G^i$ for any
$i < d$. That first implies that $\rho \otimes \psi_P$ has no trivial
subrepresentations, so $\dim H^0_{\loc}(M \otimes \calL_P) = 
\dim H^0_{\loc}(M^\dual 
\otimes \calL_{-P}) = 0$;
by Poincar\'e duality, we also have
$\dim H^1_{\loc}(M \otimes \calL_P) = 0$.
It next implies that
$\Swan(\rho \otimes \psi_P) = d \rank(M)$, so
by Theorem~\ref{thm:gos} we compute
\begin{align*}
\chi(M \otimes \calL_{P}) &= \chi(\AAA^1) \rank(M) - 
\Swan_\infty(M) \\ 
&= \rank(M) - d \rank (M) = (1-d) \rank(M).
\end{align*}
Since $H^0(M \otimes \calL_P)$ injects into $H^0_{\loc}(M 
\otimes \calL_P)$ by the exactness of \eqref{eq:snake}, it also vanishes,
and we conclude $\dim_K H^1(M \otimes \calL_P) = (d-1) \rank(M)$,
as desired.
\end{proof}

In particular, when $n$ is sufficiently large and not divisible by $p$, 
for any $r,s$ in the ring of
integers $\gotho'$ of a finite extension $K'$ of $K$, with $r$
not in the maximal ideal of $\gotho'$, we have that
$\dim_{K'} H^0(M \otimes \calL_{rx^n + s}) = 
\dim_{K'} H^0(M^\dual \otimes \calR_{-rx^n-s}) = 0$, while
$\dim_{K'} H^1(M \otimes \calL_{rx^n + s})$ and
$\dim_{K'} H^1(M^\dual \otimes \calL_{-rx^n - s})$ are equal to each
other and to a common value not depending on $r$ or $s$.
Hence we may apply Proposition~\ref{prop:irred} to deduce that
the Fourier transform of
$M \otimes \calL_{rx^n}$ is also a $(\sigma, \nabla)$-module over
$K \langle s \rangle^\dagger$. In particular, $M \otimes \calL_{rx^n}$
is (absolutely)
irreducible if and only if its Fourier transform is 
(absolutely) irreducible.

\section{Trace formulas}
\label{sec:trace}

The notion of a $p$-adic analytic Lefschetz trace formula for Frobenius first
appears in the work of Dwork \cite{bib:dwork}, \cite{bib:dwork2},
\cite{bib:dwork3}, \cite{bib:dwork4}, \cite{bib:dwork5}, and was refined by
Reich \cite{bib:reich}. The precise statement we need 
(Theorem~\ref{thm:lefschetz}; its statement also appeared
in \eqref{eq:fixpt2}) was
proved in the constant coefficient case by Monsky \cite{bib:mw3},
and in the general case by \'Etesse and le Stum
\cite[Th\'eor\`eme~6.3]{bib:etesse-lestum}.
In this chapter, we review the derivation of
\cite[Th\'eor\`eme~6.3]{bib:etesse-lestum}; since we have no real
addition to make, we will keep the review brief.

In this chapter, we will always take $k = \FF_q$, so that we can
take the Frobenius lift $\sigma_K$ on $K$ to be the identity map. In
that case, the Frobenius map $F$ on a $(\sigma, \nabla)$-module becomes
linear over $K$ (though not over a dagger algebra).

\subsection{Dwork operators}
\label{subsec:dworkop}

For $A$ a dagger algebra equipped with a Frobenius lift $\sigma$,
let $\sigma_*$ denote the ``restriction of scalars'' functor on the
category of $A$-modules. Following \cite[Definition~2.1]{bib:mw3},
a \emph{Dwork operator} on a finite $A$-module $M$ is an element of
$\Hom(\sigma_* M, M)$; such an operator can be identified with a map
$f: M \to M$ such that $f(a^\sigma \bv) = af(\bv)$ for $a \in A$ and
$\bv \in M$.

We say a Frobenius lift $\sigma: A \to A$ is \emph{Galois} if 
$\Aut(A/A^\sigma)$ has order $q$.
\begin{lemma} \label{lem:galois frob}
Let $A$ be a dagger algebra of MW-type equipped with a Galois Frobenius $\sigma$.
Then for any $(F, \nabla)$-module $M$ over $A$
and any $\tau \in \Aut(A/A^\sigma)$, there is a canonical
isomorphism $\tau^* M \to M$ of modules with connection.
\end{lemma}
Since $\tau$ reduces to the identity modulo $\gothm$,
this is a special case of the functoriality of rigid cohomology, but we
will go ahead and sketch the construction.
\begin{proof}
The desired isomorphism is constructed by ``parallel transport''.
If $A$ admits nonvanishing
local coordinates $t_1, \dots, t_n$, and $E_i$ denotes
the contraction of $\nabla$ with the vector field 
$\frac{\partial}{\partial t_i}$, the isomorphism from $\tau^* M 
= A \otimes_{\tau,A} M$ to $M$ will be
given by
\[
a \otimes \bv \mapsto a \sum_I (t_i^\tau - 1)^n 
\otimes \frac{E_1^{i_1}\cdots E_n^{i_n} \bv}{i_1!\cdots i_n!},
\]
once it is known that this series converges. In fact, if $M$ is free over
$A$, one can produce a basis of $M$ on which each $E_i$ acts via a matrix
each of whose entries has norm less than $|p|^{1/(p-1)}$: given a basis
on which $E_i$ acts via the matrix $N_i$, applying $F$ gives a basis on
which $E_i$ acts via the matrix
\[
\frac{dt_i^\sigma}{dt_i} \frac{t_i}{t_i^\sigma} N_i^\sigma,
\]
and the scalar on the left has norm less than 1. (Compare with the proof of
Lemma~\ref{lem:converge}.)

The parallel transport construction turns out to be independent of the
choice of coordinates (a routine calculation),
so it patches together to give an isomorphism of $\tau^* M$ with $M$
over all of $A$.
\end{proof}

Let $A$ be a dagger algebra of MW-type equipped with a Galois Frobenius.
Given a $(\sigma, \nabla)$-module $M$ over $A$, Lemma~\ref{lem:galois frob}
 gives
us a canonical action of $G = \Aut(A/A^\sigma)$ on $M$, and the invariants
of $M$ under $G$ are precisely $F(M)$ (as can be seen, for instance,
by expanding in series about any one point).
One can thus construct a canonical (twisted) one-sided inverse of Frobenius:
define the map $\psi: M \to M$ by the formula
\[
\psi(\bv) = F^{-1} \left( \sum_{\tau \in G} \bv^\tau \right).
\]
We also define $\psi: M \otimes \Omega^i_A \to M \otimes \Omega^i_A$ by the
same formula, using the action of $\tau$ on $\Omega^i_A$ given by
functoriality of the module of differentials. As in
\cite[Theorem~8.5]{bib:mw1} (the case $M = A$), we have that
$\psi$ is a Dwork operator,
$F \circ \psi$ equals multiplication by $q^n$ if $A$ has pure dimension
$n$, and $\psi$ commutes with $\nabla$.

\begin{remark}
One can also construct $\psi$ in the case where the Frobenius lift
$\sigma$ is not Galois; see for instance \cite[\S 4]{bib:etesse-lestum}.
\end{remark}

\subsection{Nuclearity of the canonical Dwork operator}

We next want to show that the canonical Dwork operator associated
to an $(F, \nabla)$-module is a ``trace class operator''. For this we
may appeal to the work of Monsky \cite[Section~2]{bib:mw3}; our reference
for $p$-adic functional analysis is 
\cite[Chapter~IV]{bib:schneider}.

A continuous linear map $f: V \to W$ between  locally convex $K$-vector spaces is
\emph{nuclear} if it factors as a composite $V \to V_1 \to W_1 \to W$ of
continuous linear maps, with $f_1: 
V_1 \to W_1$ a compact map between $K$-Banach
spaces. Then there is a trace functional on the set of nuclear maps from
a locally convex $K$-vector space $V$ to itself, satisfying the usual axioms:
\begin{itemize}
\item
for $c \in K$ and $f: V \to V$ nuclear,
$\Trace(cf) = c\Trace(f)$;
\item
for $f,g: V \to V$ nuclear,
$\Trace(f+g) = \Trace(f) + \Trace(g)$;
\item
for $f: V \to V$ nuclear and $g: V \to V$ continuous linear,
$f \circ g$ and $g \circ f$ are nuclear and
$\Trace(f \circ g) = \Trace(g \circ f)$.
\end{itemize}

The following result is a consequence of
\cite[Theorem~2.3]{bib:mw3} or \cite[Lemme~5.2]{bib:etesse-lestum};
we omit its proof. The reader who wishes to derive the result by hand
should start with the case $M = A = W_n$, then reduce to this case
by an analogue of the Hilbert syzygy theorem for dagger algebras.
\begin{prop}
Let $A$ be a dagger algebra of MW-type equipped with a Frobenius lift $\sigma$,
and let $M$ be a finite $A$-module. Then for each Dwork operator
$\Theta$ on $M$ (with respect to $\sigma$), $\Theta$ is nuclear,
and $\Trace(\Theta)$ is an additive function of $\Theta$.
\end{prop}

\subsection{The Lefschetz trace formula for Frobenius}

We now prove the Lefschetz trace formula \eqref{eq:fixpt2}. Again,
we follow \cite{bib:mw3}, this time imitating the ``removal of points''
trick as in \cite[Section~4]{bib:mw3}.

The following is essentially (the Galois case of)
\cite[Lemme~5.3]{bib:etesse-lestum}.
\begin{lemma} \label{lem:tracezero}
Let $A$ be a dagger algebra of MW-type equipped with a Galois Frobenius, whose
special fibre has no $\FF_q$-rational points. 
Let $M$ be a $(\sigma, \nabla)$-module over $A$, and let
$\psi$ be the canonical Dwork operator on $M$, as constructed
in Section~\ref{subsec:dworkop}. Then
$\Trace(\psi, M \otimes \Omega^i_A) = 0$ for all $i$.
\end{lemma}
\begin{proof}
We identify elements of $A$ with the continuous maps they induce on
$M \otimes \Omega^i_A$ via multiplication. Then for $a,b \in A$, we have
\begin{align*}
\Trace((a^\sigma - a)b \psi) &=
\Trace(a^\sigma b \psi) - \Trace(b a \psi) \\
&= \Trace(a^\sigma (b \psi)) - \Trace((b \psi) a^\sigma) \\
&= 0.
\end{align*}
However, by \cite[Theorem~3.3]{bib:mw3} (or an easy hand calculation),
the ideal of $A$ generated by elements of the form $a^\sigma - a$
is the unit ideal.
We conclude that the trace of $\psi$ itself on each $M \otimes \Omega^i_A$
is zero, as desired.
\end{proof}

We now deduce the trace formula as in 
\cite[Th\'eor\`eme~6.3]{bib:etesse-lestum}.
\begin{theorem} \label{thm:lefschetz}
The Lefschetz trace formula \eqref{eq:fixpt2} holds.
\end{theorem}
\begin{proof}
It suffices to check that
\begin{equation} \label{eq:trace}
\sum_{x \in X(\FF_q)} \Trace(F, \calE_x)
= \sum_i (-1)^i \Trace(F, H^i_{c,\rig}(X/K, \calE));
\end{equation}
the desired result follows by applying this assertion with $q$
replaced by each of its powers in succession.
We prove this by induction on dimension, the case $\dim(X) = 0$
being straightforward.

Both sides of \eqref{eq:trace} are additive in $X$ (the left side 
evidently, the right side by excision),
so (by the induction hypothesis) there is no loss of generality in
replacing $X$ by an open dense subset, or by an irreducible component.
In particular, we may restrict to the case where $X$ admits a finite
\'etale map to an open dense subscheme $U$ of affine $n$-space. We may also
assume that $X$ has no $\FF_q$-rational points.

Let $A$ be a dagger algebra of MW-type with special fibre $X$.
Then $A$ admits a Galois
Frobenius, by extension from the standard Frobenius on a dagger algebra
with special fibre $U$ (i.e., one extended from the standard Frobenius
on $K \langle x_1, \dots, x_n \rangle^\dagger$).
Let $M$ be a $(\sigma, \nabla)$-module over $A$ corresponding to
$\calE^\dual$, and let $\psi$ be the canonical Dwork operator on $M$.
By Lemma~\ref{lem:tracezero}, $\Trace(\psi, M \otimes \Omega^i_A) = 0$
for all $i$, and so
\[
\sum_i (-1)^i \Trace(\psi, H^i(M)) = 0
\]
as well.

Since the action of $F$ on each $H^i(M)$ is invertible
\cite[Proposition~2.1]{bib:etesse-lestum},
$\psi$ acts
on $H^i(M)$ via $q^n F^{-1}$. On the other hand, by Poincar\'e duality,
the action of $q^n F^{-1}$ on $H^i(M)$ is the transpose of the action
of $F$ on $H^{2n-i}_{c,\rig}(X/K, \calE)$, 
and in particular has the same trace.
We thus conclude that
\[
\sum_i (-1)^i \Trace(F, H^i_{c,\rig}(X/K, \calE)) = 0,
\]
which is precisely \eqref{eq:trace} because $X$ was taken to have no
$\FF_q$-rational points. The desired result follows.
\end{proof}

\section{Cohomology over finite fields}
\label{sec:arch}

With the geometric setup in place, we now introduce the archimedean
considerations that will yield our analogue of Weil II.
Much of the basic work has been carried out by Crew \cite{bib:crew1},
\cite{bib:crew2}; for the sake of the reader (and the author!) unfamiliar
with the Weil II formalism, and in keeping with our general approach,
we redo the proofs of some results from
\cite{bib:crew1} and \cite{bib:crew2} in the limited generality in 
which we need them. 

In this chapter, we again take $k = \FF_q$ and $\sigma_K$ to be the
identity map. Also, unless otherwise
specified, all curves will be smooth, geometrically irreducible, affine,
and defined over $\FF_q$.

We will always let $\iota$ denote an embedding $K^{\alg} \hookrightarrow
\CC$. As in \cite{bib:deligne}, this canard is really just a technical
convenience, but one whose removal would make the exposition substantially
more awkward.

\subsection{Weights and determinantal weights}

In this section, we introduce the notions of weights and determinantal
weights, following Crew \cite{bib:crew1}.

Suppose $q' = q^a$, and $K'$ is the smallest
unramified extension of $K$ whose residue field contains $\FF_{q'}$.
For $T: V \to V$ an endomorphism of a finite dimensional
$K'$-vector space, we say that
\begin{itemize}
\item
$T$ is \emph{$\iota$-pure of weight $w$}
if for each eigenvalue $\alpha$ of $T$,
we have $|\iota(\alpha)| = q^{(w/2)a}$;
\item
$T$ is \emph{weakly $\iota$-mixed of weight $\geq w$}
(resp.\ $\leq w$)
if for each eigenvalue $\alpha$ of $T$,
we have $|\iota(\alpha)| = q^{((w+i)/2) a}$
for some real number $i = i(\alpha) \geq 0$ (resp.\ $i \leq 0$);
\item
$T$ is \emph{strongly $\iota$-mixed} (or simply \emph{$\iota$-mixed}) 
\emph{of weight $\geq w$}
(resp.\ $\leq w$)
if for each eigenvalue $\alpha$ of $T$,
we have $|\iota(\alpha)| = q^{((w+i)/2) a}$
for some \emph{integer} $i = i(\alpha) \geq 0$ (resp.\ $i \leq 0$);
\item
$T$ is \emph{$\iota$-real} if the characteristic polynomial of $T$
has coefficients which map under $\iota$ into $\RR$. In other words,
the eigenvalues of $T: V \otimes_{\iota} \CC \to V \otimes_{\iota}
\CC$ occur in complex conjugate pairs.
\end{itemize}
If $\calE$ is an overconvergent $F$-isocrystal on a smooth
$\FF_q$-variety $X$, then we say that
$\calE$ has one of the above properties if the linear transformation 
$F_x$
on $\calE_x$ has that property for each closed point $x$ of $X$,
when we take $\FF_{q'} = \kappa(x)$.
This immediately implies that $H^0_{\rig}(X/K, \calE)$ has the same property,
since the action of Frobenius on its elements can be read off by restricting
them to any fibre of $\calE$.

We say that $\calE$ is 
\emph{$\iota$-realizable} if $\calE$ is a direct summand
of an $\iota$-real overconvergent $F$-isocrystal. Note that if $\calE$
is pure of some weight $w$, then $\calE$ is $\iota$-realizable, since
$\calE \oplus \calE^\dual(-w)$ is $\iota$-real.

\begin{remark} \label{rem:pointwise}
Beware that the definitions above correspond to what would be called
``pointwise $\iota$-pure'' and so on in \cite{bib:deligne}. In particular,
Deligne's definition of ``$\iota$-mixed'' is global and not pointwise; it
requires that $\calE$ have a filtration whose successive quotients are each
$\iota$-pure.
\end{remark}

We recall a result that
that gives us a handle on
weights in the rank one case.
The following is a
special case of a result of Tsuzuki \cite[Proposition~7.2.1]{bib:tsu4}.
(The case where $X$ is a curve is due to Crew \cite{bib:crew0}, who
showed that in that case, one may take $f$ to be finite \'etale.
This is only known in the general case modulo resolution of singularities.)
\begin{prop} \label{prop:crew}
Let $\calE$ be an overconvergent $F$-isocrystal of rank $1$ on 
a smooth $k$-variety $X$.
Then there exists a proper, dominant, generically finite \'etale morphism 
$f: Y \to X$ such that $f^* \calE$ is constant on $Y$.
\end{prop}
\begin{cor} \label{cor:const}
Let $\calE$ be an overconvergent $F$-isocrystal of rank $1$ on a 
smooth $k$-variety $X$.
Then there exists a positive integer $n$ such that $\calE^{\otimes n}$
is constant on $X$. 
\end{cor}
\begin{proof}
Choose $f: Y \to X$ as in Proposition~\ref{prop:crew}, but also
Galois (i.e., take normal closure if needed).
Let $U$ be the open dense subset of $X$ over which $f$ is \'etale,
and put $V = f^{-1}(U)$.
Let $G$ be the Galois group of $f$ and
put $n = \deg(f) = \#G$. Then $G$ acts on $H^0_{\rig}(V/K, f^* \calE)$ via
some character $\chi$, and on $H^0_{\rig}(U/K, f^* \calE^{\otimes n})$ via
$\chi^n$. The latter is the trivial character, so any horizontal
section of $f^* \calE^{\otimes n}$ descends to $\calE^{\otimes n}$,
forcing the latter to be constant on $U$. 
Moreover, by a theorem of \'Etesse \cite[Th\'eor\`eme~4]{bib:etesse},
any horizontal
section of $\calE^{\otimes n}$ over an open dense subset of $X$
extends to $X$.
This yields the desired result.
\end{proof}
In the proof of Corollary~\ref{cor:const},
one could also get by without the theorem of \'Etesse by only
proving the claim over an open dense subset of $X$, as this would
suffice for our application. We decided instead to assert the
cleaner statement.
\begin{cor}
An overconvergent $F$-isocrystal of rank $1$ on a smooth irreducible
$k$-variety is $\iota$-pure of some weight.
\end{cor}

If $\calE$ is absolutely 
irreducible of rank $d$, we define the \emph{$\iota$-determinantal
weight} of $\calE$ as $1/d$ times the $\iota$-weight of $\wedge^d \calE$
(which is unambiguous because $\wedge^d \calE$ has rank 1).
For general
$\calE$, we define the $\iota$-determinantal weights of $\calE$ to be the
$\iota$-determinantal weights of the absolutely irreducible (Jordan-H\"older)
constituents of $\calE$. (That is, extend $\FF_q, k, K$ as needed,
perform the decomposition, and define the determinantal weights, being
careful to normalize properly.)
If all of these are equal to $\alpha$, we say
$\calE$ is \emph{purely of $\iota$-determinantal weight $\alpha$}.

Since determinantal weights are defined in terms of irreducible components,
they are not \emph{a priori} well-behaved with respect to tensor products.
Our next goal is to show that they actually behave like the valuations
of eigenvalues with respect to tensoring.

\subsection{Global monodromy and determinantal weights}

So far, we have only succeeded in imposing archimedean constraints
on isocrystals of rank one. To go further, we need to use the monodromy
formalism of \cite[Chapter~1]{bib:deligne}, which we develop following
(and abbreviating) \cite[Section~4]{bib:crew1}; this formalism will allow us
to exploit results about algebraic groups to constrain first determinantal
weights and then weights. More specifically, the formalism of Tannakian
categories, which may at first seem like a contentless abstraction,
serves the vital function of giving us a mechanism for simultaneously
controlling the action of Frobenius on all tensor powers of a single
$(\sigma, \nabla)$-module. The only non-formal input needed is the 
definition of the determinantal weights, which relies on 
Proposition~\ref{prop:crew}.

Let $X$ be a smooth irreducible 
affine $\FF_q$-variety containing a rational point $x$,
and let $A$ be a dagger algebra of MW-type with special fibre $X$.
Let $\calE$ be an overconvergent $F$-isocrystal on $X$,
corresponding to a $(\sigma, \nabla)$-module $M$ over $A$;
assume for simplicity that $\calE$ is absolutely semisimple.
Note that this means $M$ is semisimple
as a module with connection: after any base extension, the sum of all
irreducible submodules with connection is an $F$-stable submodule, so
must be all of $M$. 

Let $C_M$ be the category of finite locally free $A$-modules with connection
which are isomorphic, as modules with connection, to 
subquotients of $M^a \otimes (M^\dual)^b$ for some nonnegative integers $a$
and $b$. (We crucially \emph{do not} assume they are stable under Frobenius.)
Note that the formation of the
space of morphisms between a pair of objects commutes with base extension,
because it is a finite dimensional vector space over $K$ (i.e., it is
determined by $K$-linear conditions). In particular, the irreducible
Jordan-H\"older constituents of any object in $C_M$ are absolutely 
irreducible.
Let $G_M$ be the subgroup of $\GL(M_x)$, in the category of affine algebraic groups
over $K$, of elements which commute with morphisms in $C_M$.

By construction, $G_M$ acts on $N_x$ for each $N \in C_M$;
we can read off some properties of this action from the construction.
For one, the induced action of $G_M$ on $(M^a \otimes (M^\dual)^b)_x$
must respect the decomposition of $M^a \otimes (M^\dual)^b$ into irreducibles.
In particular, $G_M$ must act trivially on any copies of the trivial
representation in $M^a \otimes (M^\dual)^b$.

It turns out that the relationship between $C_M$ and $G_M$ is much stronger
than is made evident by the trivial observations above; one has the
following result, whose proof is also essentially formal but somehow
much subtler.
\begin{prop}[Tannaka duality] \label{prop:tannaka}
The functor from $C_M$ to the category of finite dimensional
representations of $G_M$ taking $N \in C_M$ to $N_x$ is an equivalence of
categories.
\end{prop}
\begin{proof}
We only sketch the proof, since the result is well-known (see below).
In the ind-category of $C_M$ (whose elements are direct limits 
of elements of $C_M$), one can construct
an object $B$ such that $\Hom(N, B)$ is canonically isomorphic to (the underlying
set of) $N_x$. 
From this characterization, it follows that
$B_x$ (which is a direct limit of finite dimensional
representations of $C_M$) has a natural $K$-algebra structure, as well as
a coassociative, counital comultiplication. 
Thus it is the coordinate ring
of an algebraic group $G'$, whose points over any field $K'/K$
are the set of grouplike elements of $(B_x) \otimes_K K'$. (Remember that
$B_x$ is ind-finite, so tensor product commutes with its formation.)
Note that points of $G'$ act faithfully on fibres of elements of $C_M$;
that is, there is a natural functor from $C_M$ to the category of 
$B_x$-comodules
of finite type. It is not hard to show that this functor is an equivalence
of categories.

To conclude, we need only show that $G' = G_M$. This is also not difficult:
on one hand, since (points of) $G'$ acts faithfully on $C_M$, we must have
$G' \subseteq G_M$. On the other hand, $G_M$ acts on the $K$-algebra $B_x$,
but the automorphisms of the latter are given precisely by $G'$. (Namely,
the image of the identity element under an automorphism gives an element of $G'$
which induces that automorphism.) Thus $G' = G_M$, and we have the desired result.
\end{proof}

For the reader already familiar with Tannakian categories, it is possible that
 the above ``explication'' may have made things more confusing, rather than
less; consequently it is worth explaining what just happened in proper 
Tannakian language.
The category $C_M$ is a Tannakian category over $K$,
and the functor $N \to N_x$ from $C_M$ to $K$-vector spaces
is a fibre functor (so $C_M$ is actually neutral). The group $G_M$ is the
automorphism group of this fibre functor, so Proposition~\ref{prop:tannaka} is 
simply the fact that $C_M$ is equivalent to the category of finite dimensional
representations of its fibre functor, which is just an instance of Tannaka duality
\`a la \cite[Th\'eor\`eme~II.4.1.1]{bib:saavedra}.
Moreover, the proof sketch of Proposition~\ref{prop:tannaka}
is mainly a transcription of the proof of the key
intermediate result \cite[Th\'eor\`eme~II.2.3.2]{bib:saavedra}
(see also \cite[Th\'eor\`eme~II.2.6.1]{bib:saavedra}, or
\cite[Section~2.5]{bib:springer} for a discussion in more concrete
terminology).

Although Proposition~\ref{prop:tannaka} has the powerful consequence that
\emph{every} finite dimensional representation of $G_M$ occurs in the fibre
of some element of $C_M$, we will not use that fact. Instead, we need only the
much weaker result that for any irreducible $N \in C_M$, the induced action of $G_M$ on $N_x$
is irreducible, and in fact absolutely irreducible (since any submodule of $N$
over an extension of $K$ descends).

The power of the Tannakian construction is that it allows us to bring
facts about algebraic groups to bear against $M$, as follows. (The
result does not require characteristic zero, but we include the
hypothesis for simplicity.)
\begin{prop} \label{prop:outer}
Let $V$ be a finite dimensional vector space over a field $K$
of characteristic zero,
and let $G$ be an algebraic subgroup of $\GL(V)$ which acts absolutely
semisimply on $V$. Let
$N$ and $Z$ be the normalizer and centralizer, respectively,
of $G$ in $\GL(V)$. Then $N/GZ$ is finite; in other words, the group of
outer automorphisms of $G$ induced by $N$ is finite.
\end{prop}
\begin{proof}
We may as well assume that $K$ is algebraically closed, and that
the action on $G$ is irreducible; in particular, that means that by
Schur's lemma, $Z$ is the
group of scalar matrices.
Since $K$ is of characteristic zero, we can reduce to a statement about
Lie algebras as follows. 
Let $\gothgl, \gothsl, \gothg, \gothn, \gothz$ be the tangent spaces at
the identity of $\GL(V), \SL(V), G, N, Z$, respectively; we identify elements
of these spaces with linear transformations on $V$.
Then the claim is precisely that we have an equality of Lie
algebras $\gothn = \gothg + \gothz$.

To prove this, we use two facts from basic Lie theory
(see for instance \cite[Appendix~D]{bib:fh}).
\begin{enumerate}
\item[(a)] Cartan's criterion: if a Lie subalgebra $\gothg$ of $\gothgl$
(over a field of characteristic zero) satisfies $\Trace(xy) = 0$ for all
$x,y \in \gothg$, then every element of $[\gothg, \gothg]$ acts on $V$ via a
nilpotent matrix.
\item[(b)] Engel's theorem: if a Lie subalgebra $\gothg$ of $\gothgl$
has the property that every element of $\gothg$ acts on $V$ via
a nilpotent matrix, then there is a nonzero element of $V$ annihilated by
every element of $\gothg$.
\end{enumerate}
Define the trace pairing on $\gothgl$ by $x \cdot y = \Trace(xy)$; one checks that
\[
[x,y] \cdot z = y \cdot [z,x].
\]
We first observe that $\gothg \cap \gothsl$ is nondegenerate under the trace pairing
(Cartan's criterion for semisimplicity), as follows.
Let $\gothh$ be the set of $x \in \gothg$ such that $x\cdot y = 0$
for all $y \in \gothg$; then $[x,y] \in \gothh$ whenever
$x \in \gothg$ and $y \in \gothh$, which is to say
$\gothh$ is an ideal of $\gothg$. In particular, $\gothh$ is a Lie
algebra in its own right; by Cartan's criterion and Engel's theorem,
the subspace of $V$ annihilated by $[\gothh,\gothh]$ is nontrivial.
But this subspace is $\gothg$-stable because 
$[\gothh,\gothh]$ is an ideal of $\gothg$; by irreducibilty,
this forces $[\gothh,\gothh] = 0$, that is, $\gothh$ is an abelian ideal of $\gothg$.
This means the elements of $\gothh$ have joint eigenspaces, which
are again $\gothg$-stable; this can only happen if there is only one eigenspace,
that is, if $\gothh$ consists of scalar matrices.
Hence $\gothh \subseteq \gothz$ and
$\gothg \cap \gothsl$ is indeed nondegenerate under the trace pairing.

Now let $\gothg^\perp$ denote the orthogonal complement of $\gothg$
under the trace pairing. Since $\gothsl$ is itself nondegenerate under
the trace pairing (same proof as above), so is $\gothgl$; hence $\dim \gothg
+ \dim \gothg^{\perp} = \dim(V)^2$.
By $\gothg \cap \gothg^\perp \subseteq \gothz$ from above, plus the fact
that $\gothz \cap \gothz^{\perp} = 0$, we have
$\gothgl = \gothg \oplus \gothg^\perp$ as vector spaces.
It now suffices to show that if $n \in \gothn \cap \gothg^\perp$,
then $n \in \gothz$. 

By the definition of $\gothn$, for any $x \in \gothg$
we have $[n,x] \in \gothg$. But if $y \in \gothg$, we then have
$[n,x]\cdot y = n\cdot [x,y] = 0$ because $n \in \gothg^{\perp}$, 
whereas if $y \in \gothg^\perp$,
then $[n,x]\cdot y = 0$. 
Since $\gothgl$ is nondegenerate under the pairing, 
$[n,x] = 0$ and so
$n$ belongs to the centralizer of $\gothg$, which by Schur's lemma
is $\gothz$. This yields the desired result.
\end{proof}
The condition on $G$ implies that the group is reductive.
It turns out $G_M^0$ (the connected component of the identity)
is not just reductive, but semisimple
\cite[Corollary~4.10]{bib:crew1}, but we will not need this. (The 
semisimplicity is not formal, as it relies on Proposition~\ref{prop:crew}.)

For $n$ a positive integer, let $W_n$ be the semidirect product of $G_M$
by the group generated by $F^n$, and let $\deg: W_n \to \ZZ$ be the
projection onto the group generated by $F^n$ followed by the homomorphism
taking $F^{mn}$ to $mn$.
Then for some $n$, the action of $F^n$ on $G_M$ is an inner
automorphism by Proposition~\ref{prop:outer},
 so $W_n$ splits as the product of
$G_M$ by $\ZZ$.
In other words, some power
of $F$ respects the Jordan-H\"older decompositions \emph{simultaneously}
of all of the $M^a \otimes (M^\dual)^b$. (We reiterate that this simultaneity
is the principal contribution of the Tannakian point of view.)

This splitting, together with the finiteness of $\det(\rho)$, give us
the following characterization of determinantal weights, which will
conclude our consideration of $G_M$. Note that everything so far has
been purely formal; here we need the non-formal input of 
Proposition~\ref{prop:crew} via Corollary~\ref{cor:const}.
\begin{prop} \label{prop:mondetwt}
For $N \in C_M$, $N$ is purely of
$\iota$-determinantal weight $\alpha$ if and only if for some
(any) central element $z$ of $W_1$ of degree $m>0$, each eigenvalue
$\lambda$ of $z$ on $N_x$ satisfies $|\iota(\lambda)| = q^{m\alpha/2}$.
\end{prop}
\begin{proof}
There is no loss of generality in assuming $N$ is absolutely irreducible
(possibly after extending $K$).
Then by Schur's lemma, 
$z$ acts on $N_x$ by a scalar matrix, so $|\iota(\lambda)|$
is the same for all eigenvalues $\lambda$ of $z$ on $N_x$. We may thus
replace $N$ by its top exterior power $P$, which has rank one.

By Corollary~\ref{cor:const}, $P^{\otimes l}$
is constant for some positive integer $l$.
Thus all of $G_M$ acts trivially on $P_x^{\otimes l}$, so $z$ and $F^m$
have the same action there. On $P$, this means that the actions of $z$
and $F^m$ differ by an $l$-th root of unity, so the eigenvalue 
$\lambda$ of $z$ on $P_x$ satisfies $|\iota(\lambda)| = q^{m \alpha/2}$
if and only if $P$ is $\iota$-pure of weight $\alpha$, or equivalent is purely
of $\iota$-determinantal weight $\alpha$.
\end{proof}

\subsection{Determinantal weights and Dirichlet series}

The interpretation of determinantal weights provided by
Proposition~\ref{prop:mondetwt} immediately yields the
following result, which is precisely \cite[Proposition~5.7]{bib:crew1}.
Its $\ell$-adic analogue is \cite[Proposition~1.3.13]{bib:deligne}.
\begin{prop} \label{prop:detwt}
Let $X$ be a smooth irreducible affine $\FF_q$-variety,
and let $\calE, \calF$ be overconvergent $F$-isocrystals on $X$.
\begin{enumerate}
\item[(i)] If $f: Y \to X$ is a finite morphism, then
$\calE$ is purely of $\iota$-determinantal weight $\alpha$ if and only
if $f^* \calE$ is.
\item[(ii)] If $\calE$ and $\calF$
are purely of $\iota$-determinantal weights $\beta$ and $\gamma$,
then $\calE \otimes \calF$ 
is purely of $\iota$-determinantal weight $\beta + \gamma$.
\item[(iii)] If $n(\beta)$ is the sum of
the ranks of the constituents of $\calE$ of $\iota$-determinantal weight
$\beta$, then the $\iota$-determinantal weights of $\wedge^d \calE$ are
the numbers $\sum_{\beta} a(\beta) \beta$, for all collections
$a(\beta)$ of integers with $0 \leq a(\beta) \leq n(\beta)$ for all
$\beta$ and $\sum_{\beta} a(\beta) = d$.
\end{enumerate}
\end{prop}

%Using this result in lieu of \cite[Proposition~1.3.13]{bib:deligne},
%one can transcribe the proof of \cite[Th\'eor\`eme~1.5.1]{bib:deligne}
%to obtain a statement about weights, as follows.

The following lemma parallels \cite[Lemme~1.5.2]{bib:deligne}. It is
here that we first gain a real archimedean handle on Frobenius.
\begin{lemma} \label{lem:weakly mixed}
Let $\calE$ be an 
$\iota$-real overconvergent $F$-isocrystal on a
curve $X$,
and let $r$ be the largest of its $\iota$-determinantal weights.
Then $\calE_x$ is weakly $\iota$-mixed of weight $\leq r$.
\end{lemma}
\begin{proof}
The trace formula \eqref{eq:fixpt2} yields for each positive integer $i$
the equality
\begin{equation} \label{eq:purity}
\prod_{x \in X} \iota \det(1 - F_x t^{\deg(x)}, \calE_x^{\otimes 2i})^{-1}
= \frac{\iota \det(1 - F t, H^1_{c,\rig}(X/K, \calE^{\otimes 2i}))}
{\iota \det(1 - F t, H^2_{c,\rig}(X/K, \calE^{\otimes 2i}))}.
\end{equation}
In this formula, the factor
$\iota \det(1 - F_x t^{\deg(x)}, \calE_x^{\otimes 2i})^{-1}$
is a power series with nonnegative
real coefficients and constant term 1, because
the coefficient of $t^{n \deg(x)} $ is
\[
\Trace(F_x^n, \calE_x^{\otimes 2i}) = \Trace(F_x^n, \calE_x^{\otimes i})^2
\geq 0.
\]
Hence for any eigenvalue $\alpha$ of $F_x$ on $\calE_x$,
the left side of \eqref{eq:purity} (viewed as a power series about $t=0$)
has radius of convergence less than or equal to
$|\iota(\alpha^{-2i/\deg(x)})|$.

On the other hand, the radius of convergence of the right side of
\eqref{eq:purity} equals the smallest 
$\iota$-norm of any of its poles,
which is at least the smallest $\iota$-norm of any of the
eigenvalues of $F^{-1}$ on $H^2_{c,\rig}(X/K, \calE^{\otimes 2i})$.
By Poincar\'e duality, these eigenvalues are also eigenvalues of $F$ on
$H^0_{\rig}(X/K, \calE^{\otimes -2i})(1)$.
By Proposition~\ref{prop:detwt}(ii), the $\iota$-determinantal weights
of $\calE^{\otimes -2i}$ are at least $-2ir$; hence
$H^0_{\rig}(X/K, \calE^{\otimes -2i})$ is weakly $\iota$-mixed of weight
$\geq -2ir$.

In other words, the eigenvalues of $F$ on
$H^2_{c,\rig}(X/K, \calE^{\otimes 2i})$ have $\iota$-norm at least
$q^{-1-ir}$. Therefore $|\iota(\alpha^{-2i/\deg(x)})| \geq q^{-1-ir}$,
so
\[
|\iota(\alpha)| \leq q^{\deg(x)(ir+1)/(2i)}
\]
and the result follows by taking limits as $i \to \infty$.
\end{proof}

We now deduce the desired conclusion about weights.
\begin{theorem} \label{thm:realpure}
The constituents of an $\iota$-real overconvergent $F$-isocrystal
on a curve $X$ are all $\iota$-pure.
In particular, any irreducible $\iota$-realizable overconvergent $F$-isocrystal
$X$ is $\iota$-pure of some weight.
\end{theorem}
\begin{proof}
Let $\calE$ be an $\iota$-real overconvergent $F$-isocrystal on $X$.
For any $\beta$, let $\calE_\beta$ be the sum of the constituents of
$\calE$ of $\iota$-determinantal weight $\beta$. For each
$\gamma > \beta$, let $n(\gamma)$ be the sum of the ranks of the
constituents of $\calE$ of $\iota$-determinantal weight $\gamma$,
and put $N = \sum_{\gamma>\beta} n(\gamma)$. Then
the highest $\iota$-determinantal weight of $\wedge^{N+1} \calE$
is $\beta + \sum_{\gamma>\beta} n(\gamma) \gamma$, so
by Lemma~\ref{lem:weakly mixed},
$\wedge^{N+1} \calE$ is weakly $\iota$-mixed of weight
$\leq \beta + \sum_{\gamma>\beta} n(\gamma) \gamma$.

On the other hand, for each closed point $x$ and each
eigenvalue $\lambda$ of $F_x$ on $(\calE_\beta)_x$,
one of the eigenvalues of $\wedge^{N+1} \calE$ is $\lambda$ times
the determinants of the constituents of $\calE$
of $\iota$-determinantal weights greater than $\beta$. Therefore
$|\iota(\lambda)| \leq q^{\beta \deg(x)/2}$, and so
$\calE_\beta$ is weakly $\iota$-mixed of weight $\leq \beta$.

By the same reasoning applied to $\calE^\dual$,
$(\calE_\beta)^\dual$ is weakly $\iota$-mixed of weight
$\leq -\beta$. Therefore $\calE_\beta$ is $\iota$-pure of weight
$\beta$, as desired.
\end{proof}

\begin{remark}
\label{rem:realpure}
In the $\ell$-adic context, one can deduce the same result for
more general $X$ by restricting to a suitable curve; this amounts to an
application of Bertini's theorem. (See for instance 
\cite[Theorem~I.4.3]{bib:kw}.)
As noted in Section~\ref{subsec:looseends}, 
we expect a similar result to hold in
the $p$-adic setting, but its proof will be a bit more technical; in its
absence, we will have to be a bit careful in order to work around it.
\end{remark}

\subsection{Local monodromy}

We now give
what Katz dubs the ``weight drop lemma'' \cite[Lemme~1.8.1]{bib:deligne}
and explain its consequence in local monodromy.

\begin{lemma} \label{lem:weightdrop}
Let $\calE$ be an overconvergent $F$-isocrystal on a curve $X$
which is $\iota$-pure of weight $w$. Then $H^1_{c,\rig}(X/K, \calE)$ is
weakly $\iota$-mixed of weight $\leq w+2$.
\end{lemma}
\begin{proof}
Applying $\iota$ to \eqref{eq:fixpt2}, we obtain
\begin{equation} \label{eq:fixpt3}
\prod_{x \in X} \iota \det(1 - F_x t^{\deg(x)}, \calE_x)^{-1}
= \prod_i \iota \det(1 - F t, H^i_{c,\rig}(X/K, \calE))^{(-1)^{i+1}}.
\end{equation}
Write the left side of \eqref{eq:fixpt3} as $\prod_{x \in X}
\prod_{\alpha} (1 - t^{\deg(x)} \iota(\alpha))^{-1}$,
where $\alpha$ runs over the eigenvalues
of $F_x$ 
on $\calE_x$. The sum $\sum_{x \in X} \sum_{\alpha} |\iota(\alpha)|
|t|^{\deg(x)}$ is dominated by $\sum_n \#X(\FF_{q^n}) q^{nw/2} |t|^n$,
which in turn is dominated by a constant times $\sum_n q^n q^{nw/2} |t|^n$.
That sum converges for $|t| < q^{-(w+2)/2}$; thus the product on the left
side of \eqref{eq:fixpt3} also converges absolutely in that range, and
so has no zeroes there. 

We now compare with the right side of \eqref{eq:fixpt3}. First of all,
$H^0_{\rig}(X/K, \calE^\dual)$ is $\iota$-pure of weight $-w$; by Poincar\'e
duality, $H^2_{c,\rig}(X/K, \calE)$ is $\iota$-pure of weight $w+2$.
Thus $\iota \det(1 - Ft, H^2_{c,\rig}(X/K, 
\calE))$ does not vanish for $|t| < q^{-(w+2)/2}$.
We conclude that
\[
\iota \det(1 - Ft, H^1_{c,\rig}(X/K, \calE))
\]
does not vanish for $|t|
< q^{-(w+2)/2}$, since otherwise the right side of \eqref{eq:fixpt3} would
have a zero in that region. Therefore
$H^1_{c,\rig}(X/K, \calE)$ is weakly $\iota$-mixed
of weight $\leq w+2$, as desired.
\end{proof}
\begin{prop} \label{prop:weakmix}
Let $\calE$ be an overconvergent $F$-isocrystal on a curve $X$
which is $\iota$-pure of weight $w$. Then $H^0_{\loc}(X/K, \calE)$
is weakly $\iota$-mixed of weight $\leq w$.
\end{prop}
\begin{proof}
Recall a piece of the exact sequence \eqref{eq:snake}, rewritten
in ``geometric'' notation:
\[
H^0_{\rig}(X/K, \calE) \to H^0_{\loc}(X/K, \calE)
\to H^1_{c,\rig}(X/K, \calE).
\]
In this sequence, if $\calE$ is $\iota$-pure of weight $w$, then
$H^0_{\rig}(X/K, \calE)$ is $\iota$-pure of weight $w$ and
$H^1_{c,\rig}(X/K, \calE)$ is weakly $\iota$-mixed of weight $\leq w+2$
by Lemma~\ref{lem:weightdrop}. Thus
$H^0_{\loc}(X/K, \calE)$ is weakly $\iota$-mixed of weight $\leq w + 2$.
By the same token, $H^0_{\loc}(X/K, \calE^{\otimes n})$ is weakly
$\iota$-mixed of weight $\leq nw +2$. Since the latter contains
$H^0_{\loc}(X/K, \calE)^{\otimes n}$, 
we see that $H^0_{\loc}(X/K, \calE)$ is weakly $\iota$-mixed of weight
$\leq w + 2/n$ for any $n$, yielding the desired result.
\end{proof}

By applying the Jacobson-Morosov formalism of \cite[1.6]{bib:deligne},
we can refine the previous proposition
in imitation of the proof of \cite[Th\'eor\`eme~1.8.4]{bib:deligne};
alternatively, one can short-circuit
the formalism to proceed directly to the needed result, as we do here.

\begin{prop} \label{prop:strongmix}
Let $\calE$ be an overconvergent $F$-isocrystal on a curve $X$
which is $\iota$-pure of weight $w$.
Then $H^0_{\loc}(X/K, \calE)$
is \emph{strongly} $\iota$-mixed of weight $\leq w$.
\end{prop}
\begin{proof}
We may replace $X$ by a finite cover without loss of generality; in
particular, by the $p$-adic local monodromy theorem
(Proposition~\ref{prop:locmono}), we may reduce to
the case where $\calE$ is unipotent at each point $x$ of
$\overline{X} \setminus X$, for $\overline{X}$ a smooth compactification
of $X$. We may also enlarge $q$ and $K$ to ensure that each $x$ is rational
over $\FF_q$.

Let $A$ be a dagger algebra of MW-type 
with special fibre $X$, let $M$ be a $(\sigma, \nabla)$-module over $A$
corresponding to $\calE$, and let $A^{\inte} \hookrightarrow 
\calR^{\inte}_x$ be
an embedding corresponding to $x$ as in
Section~\ref{subsec:cohcurve}. By the unipotence hypothesis, there
exists a basis $\bv_1, \dots, \bv_n$ of $M \otimes \calR_x$ such that
the operator $E$ defined by $\nabla\bv= E\bv\otimes (dt/t)$
acts on the basis $\bv_1, \dots, \bv_n$ via a
nilpotent matrix over $K$. 
Redefine $\sigma$ to
be the Frobenius lift on $\calR_x^{\inte}$ sending $t$ to $t^q$,
and let $F$ be the corresponding Frobenius on $M \otimes \calR_x$;
then $F$ and $E$ are related by the identity $EF = qFE$.

Let $V$ be the $K$-span of the $\bv_i$, and let $N$ be the restriction
of $E$ to $V$.
Put $W = \ker(N)$ and $W_i = W \cap \im(N^i)$;
then $F$ acts on each $W_i$. Moreover, we have an $F$-equivariant injection
\[
W_i \otimes_K W_i(-i) \to (\ker N) \otimes_K V/(\ker N^i)
\]
sending $\bv \otimes \bw$ to $\bv \otimes \bx$ for any $\bx$
such that $N^i \bx = \bw$. The right side is a quotient of
$(\ker N) \otimes V$, which is contained in the kernel of $E$
on $(M \otimes M) \otimes \calR_x$.
By Proposition~\ref{prop:weakmix}, $H^0_{c,\rig}(X/K, \calE \otimes \calE)$
is weakly $\iota$-mixed of weight $\leq 2w$. We deduce that $W_i$
is weakly $\iota$-mixed of weight $\leq w-i$.

Let $N^\dual$ be the transpose of $N$ acting on $V^\dual$ (on which
Frobenius acts via the \emph{inverse} transpose of its action on $V$), and put
$W^\dual_i = \ker(N^\dual) \cap \im (N^\dual)^i$, which as above is
weakly $\iota$-mixed of weight $\leq -w-i$.
Then we have an
$F$-equivariant pairing
\[
W_i/W_{i+1} \otimes W^\dual_i(-i) \to W_i/W_{i+1} \otimes (\ker (N^\dual)^{i+1})
/ (\ker (N^\dual)^i) \to K,
\]
where the first arrow is an isomorphism defined as above,
and the second is a perfect pairing
induced by the canonical pairing $V \otimes V^\dual \to K$.
The perfectness of the resulting pairing means that
$W_i/W_{i+1}$ is weakly $\iota$-mixed of weight $\geq w-i$,
and hence $\iota$-pure of weight $w-i$.

Since $W$ has an $F$-stable filtration whose successive quotients are each
$\iota$-pure of weight $w-i$ for some nonnegative integer $i$, we conclude
that $W \cong H^0_{\loc,x}(X/K, \calE)$ is strongly $\iota$-mixed of weight
$\leq w$, as desired.
\end{proof}

As in the $\ell$-adic situation, Proposition~\ref{prop:strongmix}
can be viewed as affirming
a form of the weight-monodromy conjecture in equal characteristics;
we omit further details.

\subsection{Weil II on the affine line} \label{subsec:pur}

We can now deduce our main result on the affine line. 
We first use
the Fourier transform to derive a key special case, in which one gets
a stronger conclusion.

\begin{theorem} \label{thm:fourpur}
Let $M$ be an absolutely irreducible $(\sigma, \nabla)$-module
over $K \langle x \rangle^\dagger$
which is $\iota$-pure of weight $w$. Suppose there is a
nonnegative integer $d$ such that for any
finite extension $K'$ of $K$ and any $a$ in the ring of integers of
$K'$, if we put $M' = M \otimes_K K'$, we then have
$\dim_{K'} H^0_{\loc}(M' \otimes \calL_{ax}) = 
\dim_{K'} H^1_{\loc}(M' \otimes \calL_{ax}) = 0$
and $\dim_{K'} H^1(M' \otimes \calL_{ax}) = d$.
Then $H^1(M)$ and $H^1_c(M)$ are $\iota$-pure of weight $w+1$.
\end{theorem}
\begin{proof}
Without loss of generality, assume that $K$ contains $\pi$ such that
$\pi^{p-1} = -p$, and that $w=0$, so that $M$ and
$M^\dual$ have complex conjugate trace functions; then the same is
true of $M' \otimes \calL_{ax}$ and $(M')^\dual \otimes \calL_{-ax}$. 
By Poincar\'e duality, we have a
perfect pairing for each $a$:
\[
H^1_c(M' \otimes \calL_{ax}) \times H^1((M')^\dual \otimes \calL_{-ax})
\to H^2_c(K' \langle x\rangle^\dagger) \cong K'(-1).
\]
Given the assumption $\dim_{K'} H^0_{\loc}(M' \otimes \calL_{ax}) =
\dim_{K'} H^1_{\loc}(M' \otimes \calL_{ax}) = 0$, 
the ``forget supports'' map
$H^1_c(M' \otimes \calL_{ax}) \to H^1(M' \otimes \calL_{ax})$ in
\eqref{eq:snake} must be an
isomorphism. We thus have an $F$-equivariant perfect pairing
\begin{equation} \label{eq:pair}
H^1(M' \otimes \calL_{ax}) \times H^1((M')^\dual \otimes \calL_{-ax}) \to
K'(-1).
\end{equation}
In particular, 
$\dim_{K'} H^1((M')^\dual \otimes \calL_{-ax}) = d$ by duality.

By Proposition~\ref{prop:free}, $\widehat{M}$ and $\widehat{M^\dual}$ are 
$(\sigma, \nabla)$-modules over $K \langle x \rangle^\dagger$. Moreover,
$\widehat{M}$ and the pullback of $\widehat{M^\dual}$ by the map $x
\to -x$ have pointwise complex conjugate trace functions, so their direct
sum is $\iota$-real. By Proposition~\ref{prop:irred}, each of the two
is absolutely irreducible, and hence by Theorem~\ref{thm:realpure}
is $\iota$-pure of some weight $j$, necessarily the same for both.

In the pairing \eqref{eq:pair},
each factor on the left is $\iota$-pure of
weight $j$, and the object on the right is $\iota$-pure of
weight 2. This is only possible if $j+j=2$, so $j=1$. Thus
$\widehat{M}$ is $\iota$-pure of weight 1, as then is any
fibre, including $H^1(M) \cong H^1_c(M)$.
\end{proof}

By ``degenerating'' this purity result, we get the desired statement
on $\AAA^1$.
\begin{theorem} \label{thm:main}
Let $\calE$ be an overconvergent $F$-isocrystal on $\AAA^1$ which is
$\iota$-realizable and $\iota$-mixed of weight $\geq w$.
Then $H^1_{\rig}(\AAA^1/K, \calE)$ is $\iota$-mixed of
weight $\geq w+1$.
\end{theorem}
\begin{proof}
Let $M$ be a $(\sigma, \nabla)$-module over $K \langle x \rangle^\dagger$
corresponding to $\calE$.
There is no loss of generality in enlarging $k$, or in assuming that $M$
is absolutely irreducible; in particular, by Theorem~\ref{thm:realpure},
$M$ is $\iota$-pure of some weight, which we take to be $w$.

Choose an integer $N$ satisfying the conclusion
of Proposition~\ref{prop:swan} for $M$ and for $M^\dual$.
Choose $n > N$ not divisible by
$p$, let 
$f: K \langle x \rangle^\dagger \to K \langle s,x \rangle^\dagger$ and
$g: K \langle s \rangle^\dagger \to K \langle s,x \rangle^\dagger$
be the canonical embeddings, and define a $(\sigma,
\nabla)$-module on $K \langle s,x \rangle^\dagger$ by
\[
Q = f^* M^\dual \otimes_{K \langle s,x \rangle^{\dagger}} \calL_{sx^n}.
\]
(Geometrically, this corresponds to pulling back $M^\dual$ from $\AAA^1$
to $\AAA^2$ and twisting by a certain line bundle, as in the Fourier
transform.)
By Theorem~\ref{thm:pushfwd},
there exists a localization $A$ of $K
\langle s,s^{-1} \rangle^\dagger$ over which $R^1 g_! Q_A$
and $R^1 g_* Q_A$ are $(\sigma, \nabla)$-modules, where $Q_A = Q \otimes
A \langle x \rangle^\dagger$.
By Theorem~\ref{thm:degen}, we have an $F$-equivariant injection
\[
H^1_{c,\rig}(M^\dual) \hookrightarrow H^0_{\loc}(R^1 g_! Q_A).
\]

By Proposition~\ref{prop:swan}
and the choice of $N$, for
$K'$ a finite extension of $K$ and $a,c$
integers in $K'$ with $a$ not reducing to zero in the residue field,
$\dim_{K'} H^0_{\loc}((M')^\dual \otimes \calL_{ax^n+cx}) = 
\dim_{K'} H^1_{\loc}((M')^\dual \otimes \calL_{ax^n + cx}) = 0$ and
the $K'$-dimension of 
$H^1((M')^\dual \otimes \calL_{ax^n + cx})$ does not
depend on $c$.
By Theorem~\ref{thm:fourpur}, $H^1_c(M^\dual \otimes \calL_{ax^n})$ is
$\iota$-pure of weight $-w+1$; in other words,
$R^1 g_! Q_A$ is $\iota$-pure of weight $-w+1$.
By Proposition~\ref{prop:strongmix},
$H^0_{\loc}(R^1 g_! Q_A)$ is $\iota$-mixed of weight $\leq
-w+1$, so $H^1_{c,\rig}(M^\dual)$ is as well.
By Poincar\'e duality, $H^1_{\rig}(M)$ is $\iota$-mixed of
weight $\geq w+1$, as desired.
\end{proof}

\subsection{Rigid Weil II and the Weil conjectures}

To apply our results to arbitrary smooth
varieties, we employ the formalism of rigid cohomology. In so doing,
we recover the Riemann hypothesis component of the Weil conjectures.

As in \cite{bib:me8}, we use the following geometric lemma 
proved in \cite{bib:etale2} (and in the case of $k$ infinite perfect
in \cite{bib:etale}). This allows to reduce consideration of a
complicated isocrystal on a complicated variety to a more complicated
isocrystal on a less complicated variety, namely affine space.
Note that we already used the one-dimensional case of this result
once, in the proof of Theorem~\ref{thm:gos}.
\begin{prop} \label{prop:belyi}
Let $X$ be a smooth $k$-variety of dimension $n$, for $k$ a field
of characteristic $p>0$, and let $S$ be a zero-dimensional closed
subscheme of $X$. 
Then $X$ contains an open dense affine subvariety containing $S$
and admitting a finite \'etale morphism to affine $n$-space.
\end{prop}

\begin{theorem}[Rigid Weil II over a point] \label{thm:main1}
Let $X$ be a variety (separated scheme of finite type) over $\FF_q$,
and let $\calE$ be an $\iota$-realizable
overconvergent $F$-isocrystal on $X$.
\begin{enumerate}
\item[(a)]
If $\calE$ is $\iota$-mixed of weight $\leq w$, then for each $i$,
$H^i_{c,\rig}(X/K, \calE)$ is $\iota$-mixed of weight $\leq w+i$.
\item[(b)]
If $X$ is smooth and $\calE$ is $\iota$-mixed of weight $\geq w$,
then for each $i$,
$H^i_{\rig}(X/K, \calE)$ is $\iota$-mixed of weight $\geq w+i$.
\end{enumerate}
\end{theorem}
\begin{proof}
We prove the result (for all $q$) by induction primarily on $n = \dim X$
and secondarily on $\rank \calE$.
Before proceeding to the main argument, we give a number
of preliminary reductions.

Note that (b) follows from (a) by Poincar\'e duality. On the other hand,
using the excision
exact sequence \eqref{eq:excis2}
and the induction hypothesis, we may assume in (a) that $X$ is affine and
smooth of pure dimension $n$. By Poincar\'e duality again, we may now
reduce to proving just (b) for $X$, or for any one open dense subset of $X$.

There is no loss of generality in enlarging $K$, so we
assume $w=0$ by twisting as necessary. 
By Proposition~\ref{prop:belyi}, 
$X$ admits an open dense affine subscheme $U$ which in turn
admits a finite \'etale morphism $f: U \to \AAA^n$.
As noted earlier, we may replace $X$ by $U$ by excision;
since $H^i_{\rig}(U/K, \calE) \cong H^i_{\rig}(\AAA^n/K, f^* \calE)$, we may
in turn reduce to the case $X = \AAA^n$.

Finally, note that we may assume $\calE$ is irreducible:
if $0 \to \calE_1 \to \calE \to \calE_2 \to 0$ is a short exact
sequence of overconvergent $F$-isocrystals on $X$, then proving
(b) for $\calE_1$ and $\calE_2$ implies (b) for $\calE$ by the evident
long exact sequence in homology.

With these reductions in hand, we proceed to the main argument.
Choose a decomposition $\AAA^n \cong \AAA^1 \times \AAA^{n-1}$ and let
$f: \AAA^n \to \AAA^{n-1}$ be the associated projection. By
Theorem~\ref{thm:pushfwd}, there is an open dense subset $W$ of $\AAA^{n-1}$
on which the kernel $\calF_0$ and cokernel $\calF_1$ of the vertical
connection $\nabla_v$ on $\calE$ are overconvergent $F$-isocrystals
(and similarly for $\calE^\dual$).
By applying excision, we may reduce to the case $X = \AAA^1 \times W$.

Note that $f^* \calF_0$ is canonically isomorphic to a sub-$F$-isocrystal
of $\calE$. By the irreducibility hypothesis on $\calE$, if 
$\calF_0$ is nonzero, then $\calE = f^* \calF_0$. In this case,
$H^i_{\rig}(X/K, \calE) \cong H^i_{\rig}(W/K, \calF_0)$ and so the
desired result follows by the induction hypothesis. We may thus assume
that $\calF_0 = 0$. 

Since $\calE$ is $\iota$-realizable, we may choose an overconvergent
$F$-isocrystal $\calE'$ on $X$ such that $\calE \oplus \calE'$ is 
$\iota$-real; we may assume that $\calE'$ is semisimple, since
passing from $\calE'$ to its semisimplification does not change traces.
By shrinking $W$ if needed, we may ensure by
Theorem~\ref{thm:pushfwd} again that the kernel $\calF'_0$ and 
cokernel $\calF'_1$ of the vertical connection $\nabla_v$ on $\calE'$
are overconvergent $F$-isocrystals (and similarly for $(\calE')^\dual$).
Again, $f^* \calF'_0$ is canonically isomorphic to a sub-$F$-isocrystal
of $\calE'$; in particular, it is $\iota$-realizable, as then
is its restriction to $\{0\} \times W$. In other words, $\calF'_0$ itself
is $\iota$-realizable.

The trace formula \eqref{eq:fixpt2} shows that the trace of Frobenius
on a fibre of $\calF_1$, plus on that fibre of $\calF'_1$, minus
on that fibre of $\calF'_0$, is $\iota$-real. Since $\calF'_0$
is $\iota$-realizable, we deduce that $\calF_1 \oplus \calF'_1$
is $\iota$-realizable.

In particular, $\calF_1$ is $\iota$-realizable.
By Theorem~\ref{thm:main} applied to each fibre,
$\calF_1$ is $\iota$-mixed of weight $\geq 1$.
Thus by the induction hypothesis, $H^i_{\rig}(W/K, \calF_1)$
is $\iota$-mixed of weight $\geq i+1$ for each $i$.
By Proposition~\ref{prop:leray},
we have an $F$-equivariant exact sequence
\[
0 \to H^i_{\rig}(X/K,
f_* \calE^\dual) \to H^{i-1}_{\rig}(W/K, \calF_1);
\]
hence $H^i_{\rig}(X/K,
f_* \calE^\dual)$ is $\iota$-mixed of weight
$\geq i$ for all $i$.
That is, (b) holds on $X$, which thanks to the reductions completes
the induction.
\end{proof}

\begin{remark} \label{rem:overapoint}
We describe Theorem~\ref{thm:main1} as rigid Weil II ``over a point'' 
because Deligne's theorem also treats the relative
case. We do not attempt to give a relative theorem here for two reasons.
The more serious reason is that we do not have a category containing the
overconvergent $F$-isocrystals admitting Grothendieck's six operations,
so we are unable to even formulate a proper analogue. (As noted earlier,
the work of Caro \cite{bib:caro}
appears to be the appropriate setting for such an 
analogue.) The other reason is that
we are using a pointwise definition of $\iota$-mixedness, whereas Deligne's
theorem uses a global definition; see Remark~\ref{rem:pointwise}.
Remedying this discrepancy will require
extending Theorem~\ref{thm:realpure} to general varieties; 
see Remark~\ref{rem:realpure}.
\end{remark}

For completeness, we point out how Theorem~\ref{thm:main1} plus the
formalism of rigid cohomology imply the Weil conjectures in the following 
form.
\begin{enumerate}
\item[(a)] [Analytic continuation]
For $X$ a variety over $\FF_q$, the generating function
\[
\zeta_X(t) = \exp \left( \sum_{n=1}^\infty \#X(\FF_{q^n}) \frac{t^n}{n} 
\right)
\]
can be written as a product $\prod_{i=0}^{2 \dim X} P_i(t)^{(-1)^{i+1}}$,
where each $P(t)$ is a polynomial with integer coefficients and constant
coefficient 1.
\item[(b)] [Functional equation]
If $X$ is smooth, proper and purely of dimension $n$, then the product
representation can be chosen so that
\[
P_{2n-i}(t) = c t^j P_i(t^{-1})
\]
for some integer $j$ and some nonzero rational number $c$.
\item[(c)] [Riemann hypothesis]
If $X$ is smooth, proper and purely of dimension $n$, then the product
representation can be chosen so that each complex root of $P_i$
has reciprocal absolute value $q^{i/2}$.
\end{enumerate}
Part (a) follows from the Lefschetz trace formula \eqref{eq:trace}
and the finite dimensionality of rigid cohomology (with constant
coefficients), taking $P_i(t) = \det(1 - F t, H^i_{c,\rig}(X/K))$.
Part (b) follows from Poincar\'e duality and the fact that
$H^i_{c,\rig}(X/K) \cong H^i_{\rig}(X/K)$ when $X$ is proper.
Part (c) follows from Theorem~\ref{thm:main1}: it implies on one hand that
for any $\iota$,
$H^i_{c,\rig}(X/K)$ is $\iota$-mixed of weight $\leq i$,
and on the other hand
that $H^i_{c,\rig}(X/K)^\dual \cong H^{2n-i}_{\rig}(X/K)(n)
\cong H^{2n-i}_{c,\rig}(X/K)(n)$
is $\iota$-mixed of weight $\leq (2n-i)-2n = -i$. Hence
$H^i_{c,\rig}(X/K)$ is $\iota$-pure of weight $i$.

\subsection{The $p$-adic situation}

In closing, it is worth pointing out that one can
set up the same sort of framework using the $p$-adic valuation on $K^{\alg}$
as the weight formalism using the archimedean valuation on $\CC$.
This points up one of the benefits of having $p$-adic Weil II at hand:
one can treat both archimedean and $p$-adic valuations within the same
formalism.

We say an element $\alpha \in K^{\alg}$ has
\emph{slope} $s$ if $|\alpha| = |q^s|$; note that this means $s$ can
be any rational number, not just an integer. (The term ``slope'' derives
from the fact that the valuations of roots of polynomials are typically
computed as slopes of certain Newton polygons.) An overconvergent
$F$-isocrystal $\calE$ on an $\FF_q$-scheme $X$
is said to \emph{have slopes in the interval $[r,s]$} if for
each closed point $x$ of $X$ of degree $d$ over $\FF_q$,
the eigenvalues of $F_x$ on $\calE_x$ have slopes in the interval
$[dr, ds]$.

One cannot expect to precisely determine the slopes of the cohomology
of an overconvergent $F$-isocrystal, nor to limit them to integral
values; for instance, the cohomology
of an elliptic curve can have slopes 0 and 1 (in the ordinary case)
or 1/2 and 1/2 (in the supersingular case). The best we can do is 
limit the range of the variation as follows.
\begin{theorem} \label{thm:padic}
Let $X$ be a separated $\FF_q$-scheme of finite type of dimension $n$,
and let $\calE$ be an overconvergent $F$-isocrystal on $X$
which has slopes in the interval $[r,s]$. Then for $0 \leq i \leq 2n$,
$H^i_{c,\rig}(X/K,\calE)$ has slopes in the interval 
$[r+\max\{0,i-n\},s+\min\{i,n\}]$. If $X$ is smooth, then
$H^i_{\rig}(X/K, \calE)$ also has slopes in the interval
$[r+\max\{0,i-n\},s+\min\{i,n\}]$.
\end{theorem}
The key step in the proof of this theorem is the following innocuous-looking
lemma.
\begin{lemma} \label{lem:integral}
Let $\calE$ be an overconvergent $F$-isocrystal on an affine curve $C$
which has slopes in the interval $[0,\infty)$. Then
$H^1_{c,\rig}(C/K, \calE)$ also has slopes in the interval
$[0, \infty)$.
\end{lemma}
\begin{proof}
By the usual long exact sequence in homology, we may reduce to the case
of $\calE$ irreducible; in particular, we may assume that 
$H^0_{\rig}(C/K,\calE)
= H^0_{\rig}(C/K, \calE^\dual) = 0$ (since there is nothing to check if
$\calE$ is spanned by horizontal sections, the $H^1_{c,\rig}$ vanishing
in that case).
It suffices to show that $\Trace(F^i, H^1_{c,\rig}(C/K,\calE))$ 
has nonnegative $p$-adic valuation. By the Lefschetz trace formula
\eqref{eq:fixpt2}, we have
\[
\Trace(F, H^1_{c,\rig}(C/K,\calE))
= \Trace(F, H^2_{c,\rig}(C/K,\calE)) - \sum_{x \in C(\FF_{q})}
\Trace(F_x, \calE_x).
\]
By Poincar\'e duality,
$H^2_{c,\rig}(C/K,\calE)$ vanishes, while the trace of $F_x$ on each
$\calE_x$ has nonnegative $p$-adic valuation. This yields the desired
integrality for $i=1$; the
general result follows by repeating the argument over $\FF_{q^i}$.
\end{proof}

\begin{proof}[Proof of Theorem~\ref{thm:padic}]
We first verify the desired result for $X = \AAA^1$;
for brevity, let $M$ be a $(\sigma, \nabla)$-module over
$K \langle x \rangle^\dagger$ corresponding to $\calE$.
The case $i=0$ is straightforward, again since $H^0(M)$ embeds 
$F$-equivariantly into any fibre of $M$; ditto for $i=2$ via Poincar\'e
duality.
As for the case $i=1$, Lemma~\ref{lem:integral}
(applied after a twist) implies that $H^1_c(M)$ has slopes in the interval
$[r, \infty)$, and that
$H^1_c(M^\dual)$ has slopes in the interval $[-s, \infty)$.
By Grothendieck's specialization theorem (see, e.g.,
\cite[Proposition~5.14]{bib:me7}), $H^0_{\loc}(M)$ has slopes
in the interval $[r,s]$; by Poincar\'e duality, $H^1_{\loc}(M^\dual)$
has slopes in the interval $[-s+1, -r+1]$. Since $H^1(M^\dual)$
sits between $H^1_c(M^\dual)$ and $H^1_{\loc}(M^\dual)$ in the exact
sequence \eqref{eq:snake}, it has slopes in the interval $[-s, \infty)$.
By Poincar\'e duality, $H^1_c(M)$ has slopes in the interval
$[r, \infty) \cap (-\infty, s+1] = [r,s+1]$. Ditto for $H^1(M)$ by
the same arguments applied to $M^\dual$ plus Poincar\'e duality.

We now proceed to the general case, 
where we induct on $n = \dim X$;
we may assume $X$ is irreducible.
It suffices to consider
the case of cohomology with compact supports.
If $U$ is an open subset of $X$ and $Z = X \setminus U$,
the excision sequence \eqref{eq:excis2} traps
$H^i_{c,\rig}(X/K, \calE)$
between $H^i_{c,\rig}(U/K, \calE)$
and $H^i_{c,\rig}(Z/K, \calE)$.
Assuming the induction hypothesis and the fact that the claim holds
over $U$, the terms surrounding
$H^i_{c,\rig}(X/K, \calE)$ have slopes in the intervals
\[
[r+\max\{0,i-n\},s+\min\{i,n\}]
\qquad
\mbox{and}
\qquad
[r+\max\{0,i-\dim(Z)\},s+\min\{i,\dim(Z)\}].
\]
Since $\dim(Z) \leq n$ and $i - \dim(Z) \geq i - n$, the union of these
intervals is $[r+\max\{0,i-n\},s+\min\{i,n\}]$. In other words, to prove
the desired result over $X$, it suffices to prove it over $U$.

Now apply Proposition~\ref{prop:belyi} and excision again,
as in Theorem~\ref{thm:main1}, to reduce 
consideration to the case where
$X = \AAA^1 \times W$, $f: X \to W$ is the
canonical projection, and $R^j f_* \calE$ and $R^j f_* \calE^\dual$
are overconvergent $F$-isocrystals on $W$ for $j=0,1$.
Now we switch to considering cohomology without supports
(since we no longer need excision).
Applying the affine line case fibrewise, we see that
$R^j f_* \calE$ has slopes in the interval $[r,s+j]$ for
$j=0,1$.
By the induction hypothesis,
$H^i_{\rig}(W/K, R^0 f_* \calE)$ has slopes in the interval
$[r+\max\{0,i-n+1\},s+\min\{i,n-1\}]$ and
$H^{i-1}_{\rig}(W/K, R^1 f_* \calE)$ has slopes in the interval
$[r+\max\{0,i-n\},s+\min\{i,n\}]$.
By Proposition~\ref{prop:leray},
$H^i_{\rig}(X/K, \calE)$ thus has slopes in the interval
$[r+\max\{0,i-n\}, s+\min\{i,n\}]$,
as desired.
\end{proof}

\end{document}